\documentclass[a4paper,reqno,10pt]{amsart}
\usepackage[marginratio=1:1,height=637pt,width=395pt,tmargin=80pt]{geometry}

\usepackage{amsmath,latexsym,amssymb,amsthm}
\usepackage[nospace,noadjust]{cite}
\usepackage{stmaryrd}
\usepackage{xcolor}
\usepackage{url}

\theoremstyle{plain}
\newtheorem{thm}{Theorem}[section]
\newtheorem{proposition}[thm]{Proposition}
\newtheorem{lemma}[thm]{Lemma}
\newtheorem{corollary}[thm]{Corollary}

\theoremstyle{definition}
\newtheorem{definition}[thm]{Definition}
\newtheorem{remark}[thm]{Remark}
\newtheorem{exam}[thm]{Example}

\newcommand{\Id}{\operatorname{id}}
\newcommand{\pr}{\operatorname{pr}}
\newcommand{\g}{\mathfrak{g}}
\newcommand{\lie}{\mathfrak{l}}
\newcommand{\C}{\mathbb{C}}
\newcommand{\CEnd}{\mathfrak{gc}}
\newcommand{\olc}{\mathfrak{olc}}
\newcommand{\Homcs}{\operatorname{Hom}_{\mathrm{cs}}}
\newcommand{\Sh}{\mathrm{Sh}}
\newcommand{\ad}{\operatorname{ad}}
\newcommand{\Ad}{\operatorname{Ad}}
\newcommand{\Cur}{\operatorname{Cur}}
\newcommand{\im}{\operatorname{im}}
\newcommand{\cLeib}{\mathrm{cLeib}}
\newcommand{\omni}{\mathrm{omni}}
\newcommand{\lb}{\llbracket}
\newcommand{\rb}{\rrbracket}
\newcommand{\Md}{M_c^*}
\newcommand{\LL}{l^*\otimes1+1\otimes l}

\usepackage{hyperref}
\hypersetup{
colorlinks,
citecolor=blue,
filecolor=black,
linkcolor=blue,
urlcolor=black
}

\allowdisplaybreaks[2]
\makeatletter
\g@addto@macro\normalsize{%
\setlength\abovedisplayskip{6pt plus 3pt minus 2pt}%
\setlength\belowdisplayskip{6pt plus 3pt minus 2pt}%
\setlength\abovedisplayshortskip{2pt plus 3pt}%
\setlength\belowdisplayshortskip{5pt plus 3pt minus 2pt}}
\makeatother
\begin{document}

\title[Omni-representations of Leibniz conformal algebras]
{Omni-representations of Leibniz conformal algebras}

\author{Sami Benabdelhafidh}

\address{University of Sfax, Faculty of Sciences of Sfax,
BP 1171, 3038 Sfax, Tunisia}

\email{abdelahafidhsami41@gmail.com}



 \begin{abstract}
In this paper, we introduce the notion of an omni-representation of a
Leibniz conformal algebra $\mathfrak{g}$ on a $\mathbb{C}[\partial]$-module
$M$ as a Leibniz conformal algebra homomorphism from $\mathfrak{g}$ to the
conformal omni-Lie algebra $\mathfrak{gc}(M)\oplus M$. We show that the
right $\lambda$-action of a representation induces a Maurer--Cartan
element. When $M$ is a free $\mathbb{C}[\partial]$-module of finite rank,
it also defines, via the conformal dual $M_c^*$, a $1$-cocycle with values
in $M_c^*\otimes M\cong\mathfrak{gc}(M)$. Consequently, every such
representation gives rise to an omni-representation on
$M_c^*\otimes M$. We characterize omni-representations with values in
the graph of a conformal embedding tensor. Finally, we introduce the
omni-cohomology theory of Leibniz conformal algebras and establish the
relation between omni-cohomology groups and the usual cohomology groups.
\end{abstract}
\maketitle

\thispagestyle{empty}

\tableofcontents

\section{Introduction}\label{sec:intro}

Lie conformal algebras were introduced by Kac \cite{Kac} to provide
an axiomatic description of the singular part of the operator product
expansion of chiral fields in two-dimensional conformal field theory.
They are closely related to vertex algebras and infinite-dimensional
Lie algebras, and they also arise in the Hamiltonian formalism of
nonlinear evolution equations; see, for example, \cite{DSK}.
From an algebraic perspective, Lie conformal algebras can be viewed
as Lie algebras in a suitable pseudotensor category \cite{BDK}.
The structure theory of finite simple Lie conformal algebras and the
classification of their finite irreducible conformal modules have
been extensively studied \cite{DK,CK}.

As a natural non-skew-symmetric generalization of Lie conformal
algebras, Leibniz conformal algebras have attracted considerable
attention in recent years. They were introduced by Bakalov, Kac, and Voronov
\cite{BKV}, and their relationship with field algebras was subsequently
investigated in \cite{BK}. Leibniz pseudoalgebras provide a
multivariable analogue of Leibniz conformal algebras \cite{WuPseudo}.
Further developments include anti-Leibniz conformal algebras
\cite{Chtioui}, quadratic Leibniz conformal algebras and their central
extensions \cite{ZH}, the classification of rank-two Leibniz conformal
algebras \cite{WuRankTwo}, twisted relative Rota--Baxter operators on
Leibniz conformal algebras \cite{GuoWang}, and their homotopification and
categorification \cite{DasSahoo}.

Omni-Lie algebras were introduced by Weinstein \cite{Weinstein} as
linearizations of Courant algebroids \cite{LWX,KW}. Their main
feature is that Lie algebra structures on $V$ correspond to Dirac
structures in $\mathfrak{gl}(V)\oplus V$ that are graphs of linear maps
$V\to\mathfrak{gl}(V)$. Conformal and higher analogues were studied in
\cite{TZhang,SLZ}. Omni-Lie algebras are also useful for Leibniz
algebras \cite{Loday,LP}. Indeed, a representation of a Leibniz algebra
consists of a left and a right action. In general, it therefore cannot
be written as a single homomorphism into $\mathfrak{gl}(V)$. Liu and Sheng \cite{LS} proposed an alternative approach: they defined omni-representations as homomorphisms into the omni-Lie algebra $\mathfrak{gl}(V)\oplus V$.

The same obstruction appears in the conformal setting. A representation
of a Leibniz conformal algebra $\g$ on a $\C[\partial]$-module $M$
consists of two $\lambda$-actions, one on the left and one on the right. Unlike representations of Lie conformal algebras, it cannot usually
be described by a single homomorphism into the general conformal
algebra $\CEnd(M)$ of conformal linear maps on $M$ (recalled in
Section~\ref{sec:represent}).
In this paper, we carry the idea of Liu and Sheng over to
the conformal setting. An \emph{omni-representation} of $\g$ on $M$ is a
homomorphism $\rho:\g\to\olc(M)=\CEnd(M)\oplus M$ into the
\emph{conformal omni-Lie algebra} $\olc(M)$ of \cite{TZhang}
(Definition~\ref{def:olc}).

Our first results concern the right $\lambda$-action of a
representation. We show that it induces a Maurer--Cartan element
(Theorem~\ref{thm:MC}). Using the conformal dual $\Md$ of $M$, we also
show that, when $M$ is a free $\C[\partial]$-module of finite rank, it
defines a $1$-cocycle with values in $\Md\otimes M$, which is then
isomorphic to $\CEnd(M)$ (Proposition~\ref{prop:dual}). As a
consequence, every representation on a free $\C[\partial]$-module of
finite rank induces an omni-representation
(Theorem~\ref{thm:induced-omni-rep}).
We then study omni-representations with values in the graph of a
conformal embedding tensor. We show that such tensors are exactly the
left multiplication operators of Leibniz conformal structures
(Proposition~\ref{prop:T-is-ad}).

The cohomology of Leibniz conformal algebras was developed in
\cite{BKV,ZhangCoh}, with further applications in
\cite{Asif,DasSahoo,GuoWang}. We introduce the omni-cohomology
$H^\bullet_{\omni}(\g;\rho)$ of an omni-representation $\rho$. Let
$\rho=(T\circ\theta,\theta)$ be an omni-representation with values in a
graph. We prove that
\[
H^\bullet_{\omni}(\g;\rho)\cong H^\bullet(\g,\im(\theta))
\]
(Theorem~\ref{thm:graph-omni-cohomology}).
Our approach differs from that of \cite{LS} in the way the results are
organized.
Rather than treating the trivial, adjoint and graph cases as three
separate results, as in \cite{LS}, we interpret the trivial and adjoint
omni-representations as particular cases of omni-representations with
values in the graph of a conformal embedding tensor: the trivial case
corresponds to $T=0$ on $\C$, and the adjoint case to $T=\ad^L$. In fact, every conformal
embedding tensor is the left multiplication operator of a Leibniz
conformal algebra (Proposition~\ref{prop:T-is-ad}). Consequently, a
single result on omni-representations with values in a graph
(Theorem~\ref{thm:graph-omni-cohomology}) covers all these cases at
once.

The paper is organized as follows. In Section~\ref{sec:represent}, we
recall Leibniz conformal algebras, their representations and their
cohomology, and we show that the right $\lambda$-action is both a
Maurer--Cartan element and a $1$-cocycle. In
Section~\ref{sec:omni-represent}, we introduce omni-representations of
Leibniz conformal algebras, characterize them, and study those whose
images lie in the graph of a conformal embedding tensor. In
Section~\ref{sec:omni-cohomo}, we introduce the omni-cohomology, relate it
to the usual cohomology in the graph case, and deduce the trivial and
adjoint cases.

Throughout the paper, all vector spaces, algebras, and tensor
products are over $\C$. A $\C[\partial]$-module $M$ is a vector
space equipped with a linear operator $\partial:M\to M$.

\section{Leibniz conformal algebras, representations and the right \texorpdfstring{$\lambda$}{lambda}-action}
\label{sec:represent}

In this section, we recall Leibniz conformal algebras, their
representations and their cohomology, and fix our conventions for left
and right $\lambda$-actions. We then interpret the right
$\lambda$-action as a Maurer--Cartan element and, using the
conformal dual, as a $1$-cocycle with coefficients in $\Md\otimes M$.

\begin{definition}\label{defi:leibniz-algebra}
A \emph{(left) Leibniz algebra} is a vector space $\lie$ equipped with
a bilinear map $[\cdot ,\cdot ]:\lie\otimes\lie\to\lie$ satisfying
\begin{equation}\label{eq:leibniz-alg}
[x,[y,z]]=[[x,y],z]+[y,[x,z]],
\qquad \forall~x,y,z\in\lie.
\end{equation}
\end{definition}

\begin{definition}[{\cite{BKV,DasSahoo}}]\label{defi:leibniz-conformal}
A \emph{Leibniz conformal algebra} is a $\C[\partial]$-module $\g$
equipped with a $\C$-linear map
\[
[\cdot_\lambda\cdot]:\g\otimes\g\longrightarrow\g[\lambda],
\qquad
x\otimes y\longmapsto[x_\lambda y],
\]
satisfying
\begin{equation}\label{eq:sesq}
[(\partial x)_\lambda y]=-\lambda[x_\lambda y],
\qquad
[x_\lambda\partial y]=(\partial+\lambda)[x_\lambda y],
\end{equation}
and
\begin{equation}\label{eq:leibniz}
[x_\lambda[y_\mu z]]
=[[x_\lambda y]_{\lambda+\mu}z]
+[y_\mu[x_\lambda z]],
\qquad \forall~x,y,z\in\g.
\end{equation}
If the $\lambda$-bracket is skew-symmetric,
$[x_\lambda y]=-[y_{-\lambda-\partial}x]$,
then $\g$ is a Lie conformal algebra \cite{Kac}.

In expressions such as $[y_{-\lambda-\partial}x]$ and
$v_{-\mu-\partial}x$, the operator $\partial$ is moved to the left and
acts on the result of the operation.

\emph{Conventions on polynomiality.} As in \cite{BKV,Kac}, the
$\lambda$-bracket and all $\lambda$-actions are assumed throughout to be
polynomial in $\lambda$, that is, to take values in $\g[\lambda]$ and
$M[\lambda]$ respectively; the coefficient ring is always $\C$, and
$\lambda,\mu,\nu,\lambda_1,\lambda_2,\ldots$ are commuting indeterminates.
Formal power series (written $[[\lambda]]$) occur in exactly two places:
in cochains (Definition~\ref{defi:cohomology}), as in \cite{DasSahoo},
and in the conformal commutator \eqref{eq:cend-bracket} below when the
module is not finitely generated, in which case we speak of
\emph{formal} conformal algebras. Polynomiality is needed for every
substitution of the form $\lambda\mapsto-\lambda-\partial$ or
$\mu\mapsto\nu-\lambda$, and we point out each place where it is used;
in the formal case, we check polynomiality separately whenever a
substitution is performed.
\end{definition}

Let $(\g,[\cdot_\lambda\cdot])$ be a Leibniz conformal algebra. The
\emph{$j$-th products} $x_{(j)}y\in\g$, $j\geq0$, are defined by
\[
[x_\lambda y]=\sum_{j\geq0}\frac{\lambda^j}{j!}\,(x_{(j)}y),
\qquad\text{for all } x,y\in\g.
\]
Since $[x_\lambda y]$ is polynomial in $\lambda$, only finitely many $j$ satisfy
$x_{(j)}y\neq0$. In terms of the $j$-th products, the identities
\eqref{eq:sesq} and \eqref{eq:leibniz} are equivalent to the following:
for all $x,y,z\in\g$ and $j,m,n\geq0$,
\begin{align*}
(\partial x)_{(j)}y&=-j\,x_{(j-1)}y,\\
x_{(j)}(\partial y)&=\partial(x_{(j)}y)+j\,x_{(j-1)}y,\\
x_{(m)}(y_{(n)}z)&=\sum_{k=0}^{m}\binom{m}{k}(x_{(k)}y)_{(m+n-k)}z
+y_{(n)}(x_{(m)}z),
\end{align*}
with the convention $x_{(-1)}y=0$.

\begin{exam}\label{ex:current-general}
Let $(\lie,[\cdot,\cdot])$ be a Leibniz algebra. Then we can naturally
define a Leibniz conformal algebra structure on
$\Cur\lie=\C[\partial]\otimes\lie$ with the $\lambda$-bracket:
\[
[x_\lambda y]=[x,y]\qquad\text{for all } x,y\in\lie.
\]
$\Cur\lie$ is called the current Leibniz conformal algebra
associated with $(\lie,[\cdot,\cdot])$.
\end{exam}

\begin{exam}\label{ex:current}
Let $\lie=\operatorname{span}\{e_1,e_2\}$ be the two-dimensional Leibniz
algebra whose only nonzero bracket is $[e_2,e_2]=e_1$. Then the current
Leibniz conformal algebra $\g=\Cur\lie$ is free of rank two with the
$\lambda$-brackets:
\[
[e_{2\,\lambda}e_2]=e_1,\qquad
[e_{1\,\lambda}e_1]=[e_{1\,\lambda}e_2]=[e_{2\,\lambda}e_1]=0.
\]
$\g$ is not a Lie conformal algebra, since
$-[e_{2\,-\lambda-\partial}e_2]=-e_1\neq e_1=[e_{2\,\lambda}e_2]$.
\end{exam}

\begin{definition}\label{defi:homomorphism}
A $\C[\partial]$-linear map
$f:\g\to\g'$ between Leibniz conformal algebras is a
\emph{homomorphism} if
\begin{equation*}
f([x_\lambda y])=[f(x)_\lambda f(y)]',
\qquad \forall~x,y\in\g.
\end{equation*}
\end{definition}

\begin{definition}[{\cite{BKV,DasSahoo,ZhangCoh}}]
\label{defi:representation}
A \emph{representation} of a Leibniz conformal algebra
$(\g,[\cdot_\lambda\cdot])$ is a $\C[\partial]$-module $M$ equipped with
two $\C$-linear maps, called the left and right $\lambda$-actions,
\[
\g\otimes M\longrightarrow M[\lambda],
\qquad
x\otimes v\longmapsto x_\lambda v,
\quad\text{and}
\quad
M\otimes\g\longrightarrow M[\lambda],
\qquad
v\otimes x\longmapsto v_\lambda x,
\]
such that
\begin{align}
(\partial x)_\lambda v
&=
-\lambda(x_\lambda v),
&
x_\lambda(\partial v)
&=
(\partial+\lambda)(x_\lambda v),
\label{eq:rep-sesq-left}
\\
(\partial v)_\lambda x
&=
-\lambda(v_\lambda x),
&
v_\lambda(\partial x)
&=
(\partial+\lambda)(v_\lambda x),
\label{eq:rep-sesq-right}
\end{align}
and
\begin{align}
x_\lambda(y_\mu v)
&=
[x_\lambda y]_{\lambda+\mu}v
+
y_\mu(x_\lambda v),
\label{eq:rep1}
\\
x_\lambda(v_\mu y)
&=
(x_\lambda v)_{\lambda+\mu}y
+
v_\mu[x_\lambda y],
\label{eq:rep2}
\\
v_\lambda[x_\mu y]
&=
(v_\lambda x)_{\lambda+\mu}y
+
x_\mu(v_\lambda y).
\label{eq:rep3}
\end{align}
\end{definition}

\begin{definition}\label{defi:rep-hom}
A \emph{homomorphism} (resp.\ \emph{isomorphism}) of representations
$\phi:M\to M'$ is a $\C[\partial]$-linear map (resp.\ bijection) such
that
\[
\phi(x_\lambda v)=x_\lambda\phi(v),
\qquad
\phi(v_\lambda x)=\phi(v)_\lambda x,
\qquad\forall~x\in\g,~v\in M.
\]
\end{definition}

We now recall the cohomology of a Leibniz conformal
algebra with coefficients in a representation and the Balavoine
bracket, whose Maurer--Cartan elements are Leibniz conformal
structures \cite{Balavoine,DasSahoo,ZhangCoh}.

\begin{definition}[{\cite{BKV,DasSahoo,ZhangCoh}}]
\label{defi:cohomology}
Let $\g$ and $M$ be $\C[\partial]$-modules. For $n\geq1$, a
$\C$-linear map
\[
f:\g^{\otimes n}
\longrightarrow
M[[\lambda_1,\ldots,\lambda_{n-1}]]
\]
is called \emph{conformal sesquilinear} if
\begin{align*}
f_{\lambda_1,\ldots,\lambda_{n-1}}
(x_1,\ldots,\partial x_i,\ldots,x_n)
&=
-\lambda_i
f_{\lambda_1,\ldots,\lambda_{n-1}}
(x_1,\ldots,x_n),
\quad
1\leq i\leq n-1,
\\
f_{\lambda_1,\ldots,\lambda_{n-1}}
(x_1,\ldots,x_{n-1},\partial x_n)
&=
(\partial+\Lambda_{n-1})
f_{\lambda_1,\ldots,\lambda_{n-1}}
(x_1,\ldots,x_n),
\end{align*}
where
$\Lambda_k=\lambda_1+\cdots+\lambda_k$.
Here the parameter $\lambda_i$ is attached to the $i$-th argument
$x_i$ for $1\leq i\leq n-1$, while the last argument $x_n$ carries no
parameter; its behaviour is determined by the second sesquilinearity
condition above. We denote the space of all such maps by
\[
\Homcs\big(
\g^{\otimes n},
M[[\lambda_1,\ldots,\lambda_{n-1}]]
\big).
\]
For $n=1$, this space is
$\operatorname{Hom}_{\C[\partial]}(\g,M)$.
Let $M$ be a representation of $\g$. The space of $n$-cochains is
\[
C^n_{\cLeib}(\g,M)
=
\begin{cases}
M/\partial M, & n=0,\\[1mm]
\Homcs\big(
\g^{\otimes n},
M[[\lambda_1,\ldots,\lambda_{n-1}]]
\big), & n\geq1.
\end{cases}
\]

The coboundary map
$\delta:C^n_{\cLeib}(\g,M)\to C^{n+1}_{\cLeib}(\g,M)$ is given by
\begin{equation}\label{eq:delta0}
\delta(v+\partial M)(x)
=
-(v_\lambda x)\big|_{\lambda=0},
\qquad
v+\partial M\in M/\partial M,
\end{equation}
and, for $f\in C^n_{\cLeib}(\g,M)$ with $n\geq1$,
\begin{equation}\label{eq:delta}
\begin{aligned}
(\delta f)_{\lambda_1,\ldots,\lambda_n}(x_1,\ldots,x_{n+1})
&=\sum_{i=1}^{n}(-1)^{i+1}\,
x_{i\,\lambda_i}\,
f_{\lambda_1,\ldots,\widehat{\lambda_i},\ldots,\lambda_n}
(x_1,\ldots,\widehat{x_i},\ldots,x_{n+1})\\
&\quad+(-1)^{n+1}\,
f_{\lambda_1,\ldots,\lambda_{n-1}}(x_1,\ldots,x_n)_{\Lambda_n}\,x_{n+1}\\
&\quad+\sum_{1\leq i<j\leq n+1}(-1)^i\,
f_{\lambda_1,\ldots,\widehat{\lambda_i},\ldots,\lambda_{j-1},
\lambda_i+\lambda_j,\ldots,\lambda_n}\\
&\hspace{3em}\big(x_1,\ldots,\widehat{x_i},\ldots,x_{j-1},
[x_{i\,\lambda_i}x_j],x_{j+1},\ldots,x_{n+1}\big).
\end{aligned}
\end{equation}
When $j=n+1$, the parameter $\lambda_i+\lambda_j$ is omitted, since
the last argument of a cochain carries no conformal parameter.
Each term of \eqref{eq:delta} lies in $M[[\lambda_1,\ldots,\lambda_n]]$:
the $\lambda$-actions are polynomial in their parameter and are applied
coefficientwise (the right action being evaluated at $\Lambda_n$), and
the only substitutions are $\lambda_i\mapsto\lambda_i+\lambda_j$, which
are well defined on formal power series since each monomial receives
finitely many contributions. No substitution of the form
$\lambda\mapsto-\lambda-\partial$ occurs in the differential; such
substitutions are used only for polynomial maps, as in
\eqref{eq:endo-convention} below. One checks directly that $\delta f$ is
conformal sesquilinear, so $\delta$ maps $C^n_{\cLeib}(\g,M)$ to
$C^{n+1}_{\cLeib}(\g,M)$.
By \cite{DasSahoo,ZhangCoh}, one has $\delta^2=0$. The cohomology of
$(C^\bullet_{\cLeib}(\g,M),\delta)$ is denoted by
$H^\bullet(\g,M)$.
\end{definition}

The map \eqref{eq:delta0} is well defined since
$(\partial v)_\lambda x=-\lambda(v_\lambda x)$,
and hence its value at $\lambda=0$ vanishes. Moreover,
\[
\delta(v+\partial M)(\partial x)
=
-\big(v_\lambda\partial x\big)\big|_{\lambda=0}
=
-\partial(v_\lambda x)\big|_{\lambda=0},
\]
so that $\delta(v+\partial M)$ is a $\C[\partial]$-linear map
$\g\to M$. The identity $\delta^2=0$ in degree zero follows from
\eqref{eq:rep3} evaluated at $\lambda=0$.

We use the reduced complex of \cite{ZhangCoh,DasSahoo}, whose cochains
are not required to be skew-symmetric, and, as in \cite{LS}, the
degree-zero differential \eqref{eq:delta0} uses the right $\lambda$-action.

The space $\C$, on which $\partial$ acts by zero, with zero left and
right $\lambda$-actions is a representation of $\g$, called the
\emph{trivial representation}. We write
$H^\bullet(\g)=H^\bullet(\g,\C)$.

Another important representation is the \emph{adjoint representation}
of $\g$ on itself, whose left and right $\lambda$-actions are
\[
x_\lambda y=[x_\lambda y],
\qquad
y_\lambda x=[y_\lambda x],
\qquad \forall~x,y\in\g.
\]
The sesquilinearity conditions
\eqref{eq:rep-sesq-left}--\eqref{eq:rep-sesq-right} follow from
\eqref{eq:sesq}, and the identities \eqref{eq:rep1}--\eqref{eq:rep3}
are instances of the Leibniz conformal identity \eqref{eq:leibniz}.
We denote the corresponding cohomology by $H^\bullet(\g,\g)$.

In terms of the $j$-th
products, the degree-zero differential \eqref{eq:delta0} for the
adjoint representation reads
\[
\delta(y+\partial\g)(x)=-\,y_{(0)}x,\qquad \forall~x,y\in\g.
\]
In particular, $H^0(\g,\g)$ is
the space of classes $y+\partial\g$ such that $y_{(0)}x=0$ for all
$x\in\g$.

A permutation $\sigma\in\mathbb{S}_n$ is called an
\emph{$(i,n-i)$-shuffle} if $\sigma(1)<\cdots<\sigma(i)$ and
$\sigma(i+1)<\cdots<\sigma(n)$. If $i=0$ or $i=n$, we set $\sigma=\Id$. The
set of all $(i,n-i)$-shuffles is denoted by $\Sh(i,n-i)$.

Set $$C^\bullet_{\cLeib}(\g,\g)=\bigoplus_{k\geq1}C^k_{\cLeib}(\g,\g),$$ where
elements of $C^k_{\cLeib}(\g,\g)$ have degree $k-1$. Following
\cite{Balavoine,DasSahoo} (see \cite[Section~3.1]{YuanLiu} for the Lie
conformal case), for $f\in C^m_{\cLeib}(\g,\g)$ and
$g\in C^n_{\cLeib}(\g,\g)$, define the \emph{Balavoine bracket} on
$C^\bullet_{\cLeib}(\g,\g)$ by
\begin{equation}\label{eq:balavoine-bracket}
\lb f,g\rb=f\diamond g-(-1)^{(m-1)(n-1)}g\diamond f,
\end{equation}
where $f\diamond g\in C^{m+n-1}_{\cLeib}(\g,\g)$ is defined by
($x_1,\ldots,x_{m+n-1}\in\g$)
\begin{equation}\label{eq:partial-comp}
\begin{aligned}
&(f\diamond g)_{\lambda_1,\ldots,\lambda_{m+n-2}}(x_1,\ldots,x_{m+n-1})\\
&\quad=\sum_{i=1}^{m}\sum_{\sigma\in\Sh(i-1,n-1)}
(-1)^{(i-1)(n-1)}\operatorname{sgn}(\sigma)\,
f_{\lambda_{\sigma(1)},\ldots,\lambda_{\sigma(i-1)},\Lambda^\sigma_i,
\lambda_{i+n},\ldots,\lambda_{m+n-2}}\\
&\hspace{3em}\Big(x_{\sigma(1)},\ldots,x_{\sigma(i-1)},
g_{\lambda_{\sigma(i)},\ldots,\lambda_{\sigma(i+n-2)}}
\big(x_{\sigma(i)},\ldots,x_{\sigma(i+n-2)},x_{i+n-1}\big),\\
&\hspace{3em}\phantom{\Big(}x_{i+n},\ldots,x_{m+n-1}\Big),
\end{aligned}
\end{equation}
where $\Lambda^\sigma_i=\lambda_{\sigma(i)}+\cdots+\lambda_{\sigma(i+n-2)}+\lambda_{i+n-1}$,
if it occurs; for $i=m$ the last argument of $f$ carries no parameter, so
$\Lambda^\sigma_m$ is omitted. In particular, for $f,g\in C^2_{\cLeib}(\g,\g)$,
we have $\lb f,g\rb=f\diamond g+g\diamond f$ with
\[
(f\diamond g)_{\lambda,\mu}(x,y,z)
=f_{\lambda+\mu}(g_\lambda(x,y),z)-f_\lambda(x,g_\mu(y,z))+f_\mu(y,g_\lambda(x,z)).
\]

\begin{thm}[{\cite{Balavoine,DasSahoo}}]
\label{thm:graded-Lie}
$(C^\bullet_{\cLeib}(\g,\g),\lb\cdot,\cdot\rb)$ is a graded Lie algebra.
Moreover, $\pi\in C^2_{\cLeib}(\g,\g)$ defines a Leibniz conformal algebra
structure on $\g$ by $[x_\lambda y]=\pi_\lambda(x,y)$ if and only if
$\lb\pi,\pi\rb=0$.
\end{thm}

Indeed, $\frac12\lb\pi,\pi\rb_{\lambda,\mu}(x,y,z)
=\pi_{\lambda+\mu}(\pi_\lambda(x,y),z)-\pi_\lambda(x,\pi_\mu(y,z))
+\pi_\mu(y,\pi_\lambda(x,z))$, which vanishes if and only if
\eqref{eq:leibniz} holds. Consequently, if $\lb\pi,\pi\rb=0$, then
$d_\pi=\lb\pi,\cdot\rb$ satisfies $d_\pi^2=0$, and for
$\varpi\in C^2_{\cLeib}(\g,\g)$,
\begin{equation}\label{eq:MC-expansion}
\lb\pi+\varpi,\pi+\varpi\rb=2\Big(d_\pi\varpi+\frac12\lb\varpi,\varpi\rb\Big).
\end{equation}

We also recall some facts on conformal linear maps.
A \emph{conformal linear map} on a $\C[\partial]$-module $M$ is a
$\C$-linear map
 $A:M\longrightarrow M[\lambda], \quad v\longmapsto A_\lambda(v)$,
satisfying
$A_\lambda(\partial v)=(\partial+\lambda)A_\lambda(v)$.
The space $\operatorname{Cend}(M)$ of conformal linear maps is a
$\C[\partial]$-module via
\[
(\partial A)_\lambda(v)=-\lambda A_\lambda(v).
\]
Its conformal composition is
\[
(A_\lambda B)_\mu(v)=A_\lambda(B_{\mu-\lambda}(v)).
\]
The corresponding conformal commutator is
\begin{equation}\label{eq:cend-bracket}
[A_\lambda B]_\nu(v)
=
A_\lambda(B_{\nu-\lambda}(v))
-
B_{\nu-\lambda}(A_\lambda(v)).
\end{equation}
Equivalently,
\[
[A_\lambda B]_{\lambda+\mu}(v)=A_\lambda(B_\mu(v))-B_\mu(A_\lambda(v)).
\]
Following \cite{Kac,BKL}, we denote by $\CEnd(M)$ the space
$\operatorname{Cend}(M)$ equipped with the conformal commutator
\eqref{eq:cend-bracket}, and we call it the \emph{general conformal
algebra} of $M$. The substitution $\mu=\nu-\lambda$ in
\eqref{eq:cend-bracket} is legitimate because $A_\lambda(v)$ and
$B_\mu(v)$ are polynomial in $\lambda$ and $\mu$.

If $M$ is finitely generated, $[A_\lambda B]$ is polynomial in
$\lambda$, and $\CEnd(M)$ is a Lie conformal algebra \cite{Kac,BKL}.
In general, $[A_\lambda B]$ is only a formal power series in $\lambda$.
In this case, we call $\CEnd(M)$ a \emph{formal} Lie conformal algebra,
and linear maps are extended coefficientwise to formal power series.

Let $(M,l,r)$ be a representation of $\g$.
We use the endomorphism convention
\begin{equation}\label{eq:endo-convention}
l(x)_\lambda(v)=x_\lambda v,
\qquad
r(x)_\lambda(v)=v_{-\lambda-\partial}x,
\end{equation}
where the substitution is defined as follows. Since the right
$\lambda$-action takes values in $M[\lambda]$, we may write
$v_\nu x=\sum_{k}\nu^k g_k$, a finite sum with $g_k\in M$, and we set
\[
v_{-\lambda-\partial}x=\sum_{k}(-\lambda-\partial)^k g_k\in M[\lambda].
\]
Thus $\partial$ always acts on the coefficients $g_k$, after the right
action has been computed; every substitution of this type in the paper
is understood coefficientwise in this sense, and is well defined because
only polynomials are involved. By \eqref{eq:rep-sesq-left}, $l$ is a
$\C[\partial]$-linear map $\g\to\CEnd(M)$. For the adjoint representation, written $(\g,\ad^L,\ad^R)$, we have
\[
\ad^L(x)_\lambda(y)=[x_\lambda y],
\qquad
\ad^R(x)_\lambda(y)=[y_{-\lambda-\partial}x].
\]
More generally, for any Leibniz conformal algebra $\mathcal A$ we write
$\ad^L_{\mathcal A}:\mathcal A\to\CEnd(\mathcal A)$,
$\ad^L_{\mathcal A}(u)_\lambda(v)=[u_\lambda v]$, for its left adjoint
operator; thus $\ad^L=\ad^L_{\g}$.

For a representation $(M,l,r)$ of $\g$, it is obvious that $(M,l,0)$ is
also a representation of $\g$. Thus, we have two semidirect product
Leibniz conformal algebras $\g\ltimes_{(l,r)}M$ and $\g\ltimes_{(l,0)}M$
with the $\lambda$-brackets $[\cdot_\lambda\cdot]_{(l,r)}$ and
$[\cdot_\lambda\cdot]_{(l,0)}$ given respectively by
\begin{align}
[(x,a)_\lambda(y,b)]_{(l,r)}
&=
\big([x_\lambda y],\,x_\lambda b+a_\lambda y\big),
\label{eq:semidirect}
\\
[(x,a)_\lambda(y,b)]_{(l,0)}
&=
\big([x_\lambda y],\,x_\lambda b\big).
\label{eq:semidirect0}
\end{align}
The right $\lambda$-action induces a conformal sesquilinear map
$\bar r:(\g\oplus M)\otimes(\g\oplus M)\to(\g\oplus M)[[\lambda]]$ as
follows:
\[
\bar r_\lambda\big((x,a),(y,b)\big)=(0,a_\lambda y),
\qquad\forall~ x,y\in\g,\ a,b\in M.
\]
Let $\pi_0=[\cdot_\lambda\cdot]_{(l,0)}$. By Theorem~\ref{thm:graded-Lie},
$\big(C^\bullet_{\cLeib}(\g\oplus M,\g\oplus M),\lb\cdot,\cdot\rb,d=\lb\pi_0,\cdot\rb\big)$
is a differential graded Lie algebra. An element
$\varpi\in C^2_{\cLeib}(\g\oplus M,\g\oplus M)$ is a \emph{Maurer--Cartan
element} if $d\varpi+\frac12\lb\varpi,\varpi\rb=0$, that is, by
\eqref{eq:MC-expansion}, if $\pi_0+\varpi$ is a Leibniz conformal structure
on $\g\oplus M$.

The adjoint representation of the Leibniz conformal algebra
$\g\ltimes_{(l,0)}M$ is given by the left and right $\lambda$-actions
\[
(x,a)_\lambda(y,b)=\big([x_\lambda y],\,x_\lambda b\big),
\qquad
(y,b)_\lambda(x,a)=\big([y_\lambda x],\,y_\lambda a\big).
\]

The following theorem shows that the right $\lambda$-action of a
representation induces a Maurer--Cartan element.

\begin{thm}\label{thm:MC}
With the above notations, $\bar r$ satisfies the following
Maurer--Cartan equation on the Leibniz conformal algebra
$\g\ltimes_{(l,0)}M$:
\[
\delta\bar r-\frac12\lb\bar r,\bar r\rb=0
\quad\Longleftrightarrow\quad
d\bar r+\frac12\lb\bar r,\bar r\rb=0,
\]
where $\delta$ is the coboundary operator for the Leibniz conformal
algebra $\g\ltimes_{(l,0)}M$ with coefficients in the adjoint
representation and $d=\lb\pi_0,\cdot\rb$; the two equations are
equivalent since $d=-\delta$ on $C^2_{\cLeib}(\g\oplus M,\g\oplus M)$.
Consequently, $\pi_0+\bar r=[\cdot_\lambda\cdot]_{(l,r)}$, that is, the Leibniz conformal algebra
$\g\ltimes_{(l,r)}M$ is a deformation of the Leibniz conformal algebra
$\g\ltimes_{(l,0)}M$ via the Maurer--Cartan element $\bar r$.
\end{thm}

\begin{proof}
Let $X=(x,a)$, $Y=(y,b)$, $Z=(z,c)\in\g\ltimes_{(l,0)}M$. We first
compute both terms of the Maurer--Cartan equation directly from the
definitions, and then apply the representation identities.

\emph{Step 1: direct computation.} By \eqref{eq:delta},
\eqref{eq:semidirect0} and the definition of $\bar r$, we have
\begin{align*}
(\delta\bar r)_{\lambda,\mu}(X,Y,Z)
={}&[X_\lambda\bar r_\mu(Y,Z)]_{(l,0)}-[Y_\mu\bar r_\lambda(X,Z)]_{(l,0)}
-[\bar r_\lambda(X,Y)_{\lambda+\mu}Z]_{(l,0)}\\
&-\bar r_{\lambda+\mu}([X_\lambda Y]_{(l,0)},Z)
+\bar r_\lambda(X,[Y_\mu Z]_{(l,0)})
-\bar r_\mu(Y,[X_\lambda Z]_{(l,0)})\\
={}&\big(0,\,x_\lambda(b_\mu z)-y_\mu(a_\lambda z)
-(x_\lambda b)_{\lambda+\mu}z+a_\lambda[y_\mu z]-b_\mu[x_\lambda z]\big).
\end{align*}
Here the third term vanishes, since the $\g$-component of
$\bar r_\lambda(X,Y)$ is zero and the bracket \eqref{eq:semidirect0} has
zero $M$-component when its first argument lies in $M$.
Similarly, by the formula for $\frac12\lb\pi,\pi\rb$ following
Theorem~\ref{thm:graded-Lie}, we have
\begin{align*}
\lb\bar r,\bar r\rb_{\lambda,\mu}(X,Y,Z)
={}&2\Big(\bar r_{\lambda+\mu}\big(\bar r_\lambda(X,Y),Z\big)
-\bar r_\lambda\big(X,\bar r_\mu(Y,Z)\big)
+\bar r_\mu\big(Y,\bar r_\lambda(X,Z)\big)\Big)\\
={}&2\big(0,\,(a_\lambda y)_{\lambda+\mu}z\big).
\end{align*}
Hence
\begin{align*}
\Big(\delta\bar r-\frac12\lb\bar r,\bar r\rb\Big)_{\lambda,\mu}(X,Y,Z)
={}&\big(0,\,x_\lambda(b_\mu z)-y_\mu(a_\lambda z)-(x_\lambda b)_{\lambda+\mu}z\\
&\qquad+a_\lambda[y_\mu z]-b_\mu[x_\lambda z]-(a_\lambda y)_{\lambda+\mu}z\big).
\end{align*}

\emph{Step 2: representation identities.} Reordering the six terms
without changing their signs, the second component above is
\begin{equation*}
\underbrace{a_\lambda[y_\mu z]-(a_\lambda y)_{\lambda+\mu}z
-y_\mu(a_\lambda z)}_{=:E_1}
\;+\;
\underbrace{x_\lambda(b_\mu z)-(x_\lambda b)_{\lambda+\mu}z
-b_\mu[x_\lambda z]}_{=:E_2}.
\end{equation*}
Here $E_1$ collects the three terms involving the right action of $a$
(the term $-(a_\lambda y)_{\lambda+\mu}z$ being produced by
$-\frac12\lb\bar r,\bar r\rb$), and $E_2$ collects the three terms
involving $b$.

Identity \eqref{eq:rep3}, applied with $v=a$, with $(x,y)$ replaced by
$(y,z)$ and with the same parameters $\lambda,\mu$, reads
\[
a_\lambda[y_\mu z]=(a_\lambda y)_{\lambda+\mu}z+y_\mu(a_\lambda z),
\]
so $E_1=0$. Identity \eqref{eq:rep2}, applied with $v=b$, with $y$
replaced by $z$ and with the same parameters $\lambda,\mu$, reads
\[
x_\lambda(b_\mu z)=(x_\lambda b)_{\lambda+\mu}z+b_\mu[x_\lambda z],
\]
so $E_2=0$. Hence the second component vanishes, the first one being
zero by construction, and
$\delta\bar r-\frac12\lb\bar r,\bar r\rb=0$.

\emph{Step 3: relation with $d$.} By \eqref{eq:rep-sesq-right},
$\bar r$ is conformal sesquilinear, so
$\bar r\in C^2_{\cLeib}(\g\oplus M,\g\oplus M)$. For
$f\in C^2_{\cLeib}(\g\oplus M,\g\oplus M)$, the formula for $\diamond$ on
$C^2_{\cLeib}$ preceding Theorem~\ref{thm:graded-Lie} yields
\begin{align*}
(\pi_0\diamond f)_{\lambda,\mu}(X,Y,Z)
&=[f_\lambda(X,Y)_{\lambda+\mu}Z]_{(l,0)}-[X_\lambda f_\mu(Y,Z)]_{(l,0)}
+[Y_\mu f_\lambda(X,Z)]_{(l,0)},\\
(f\diamond\pi_0)_{\lambda,\mu}(X,Y,Z)
&=f_{\lambda+\mu}([X_\lambda Y]_{(l,0)},Z)-f_\lambda(X,[Y_\mu Z]_{(l,0)})
+f_\mu(Y,[X_\lambda Z]_{(l,0)}).
\end{align*}
Comparing with \eqref{eq:delta} for $n=2$, we obtain
\[
d f=\lb\pi_0,f\rb=\pi_0\diamond f+f\diamond\pi_0=-\delta f .
\]
Hence the equation proved above is equivalent to
$d\bar r+\frac12\lb\bar r,\bar r\rb=0$, so $\bar r$ is a Maurer--Cartan
element of the differential graded Lie algebra
$\big(C^\bullet_{\cLeib}(\g\oplus M,\g\oplus M),d,\lb\cdot,\cdot\rb\big)$.
Equivalently, $\pi_0+\bar r$ is a Leibniz conformal structure on
$\g\oplus M$; by \eqref{eq:semidirect} and \eqref{eq:semidirect0} it is
precisely $[\cdot_\lambda\cdot]_{(l,r)}$. The proof is finished.
\end{proof}

We now describe the right $\lambda$-action through the conformal dual.
Recall \cite{Kac,BKV} that the \emph{conformal dual} of a
$\C[\partial]$-module $M$ is
\[
\Md=\{f:M\to\C[\lambda]\mid f\ \C\text{-linear},\ f_\lambda(\partial v)=\lambda f_\lambda(v)\},
\quad(\partial f)_\lambda(v)=-\lambda f_\lambda(v).
\]
For convenience, we set $\langle f,v\rangle_\lambda=f_\lambda(v)$ for all
$f\in\Md$ and $v\in M$, 
so that \[\langle f,\partial v\rangle_\lambda=\lambda\langle f,v\rangle_\lambda
\quad \text{and} \quad\langle\partial f,v\rangle_\lambda=-\lambda\langle f,v\rangle_\lambda.\]
We use two hypotheses on $M$ in the rest of this section.
\begin{itemize}
\item[(H1)] $M$ is finitely generated. This hypothesis is in force until
the end of the section. It ensures that $l^*(x)_\lambda f$, defined
below, is polynomial in $\lambda$.
\item[(H2)] $M$ is free of finite rank. This stronger hypothesis is
used from the definition of $\Upsilon$ on, and in
Proposition~\ref{prop:dual}. It ensures that $\Upsilon$ is an
isomorphism.
\end{itemize}
Define
$l^*:\g\to\CEnd(\Md)$ by
\[
\langle l^*(x)_\lambda f,v\rangle_\mu=-\langle f,l(x)_\lambda v\rangle_{\mu-\lambda},\qquad
\forall ~x\in\g,\ f\in\Md,\ v\in M .
\]

Since \[\langle l^*(x)_\lambda f,\partial v\rangle_\mu
=-\langle f,(\partial+\lambda)l(x)_\lambda v\rangle_{\mu-\lambda}
=\mu\,\langle l^*(x)_\lambda f,v\rangle_\mu,\] each coefficient of $l^*(x)_\lambda f$
lies in $\Md$, and since $M$ is finitely generated, $l^*(x)_\lambda f$ is
a polynomial in $\lambda$; moreover $l^*$ is $\C[\partial]$-linear and
$l^*(x)\in\CEnd(\Md)$. By \eqref{eq:rep1},
\[
\begin{aligned}
\langle l^*(x)_\lambda(l^*(y)_\mu f)-l^*(y)_\mu(l^*(x)_\lambda f),v\rangle_\nu
&=\langle f,y_\mu(x_\lambda v)-x_\lambda(y_\mu v)\rangle_{\nu-\lambda-\mu}\\
&=\langle l^*([x_\lambda y])_{\lambda+\mu}f,v\rangle_\nu,
\end{aligned}
\]
so that $(\Md,l^*,0)$ is a representation of $\g$. The tensor product
$\Md\otimes M$ is a $\C[\partial]$-module via
$\partial=\partial\otimes1+1\otimes\partial$, and
$(\Md\otimes M,\LL,0)$ is also a representation of $\g$, where
\[
(\LL)(x)_\lambda(f\otimes u)=(l^*(x)_\lambda f)\otimes u+f\otimes x_\lambda u,
\qquad\forall f\in\Md,\ u\in M;
\]
indeed, sesquilinearity is checked directly, and in the commutator of
$(\LL)(x)_\lambda$ and $(\LL)(y)_\mu$ the mixed terms cancel, so that
\eqref{eq:rep1} follows from the corresponding identities for $l^*$ and $l$.

\begin{exam}\label{ex:dual-rep}
Let $(M,l,r)$ be a finitely generated representation of $\g$ which is
free as a $\C[\partial]$-module. Let $l^*$ and $r^*$ be the two
$\C[\partial]$-module homomorphisms from $\g$ to $\CEnd(\Md)$ given by
\[
\langle l^*(x)_\lambda f,v\rangle_\mu=-\langle f,l(x)_\lambda v\rangle_{\mu-\lambda},
\qquad
\langle r^*(x)_\lambda f,v\rangle_\mu=-\langle f,r(x)_\lambda v\rangle_{\mu-\lambda},
\]
for all $x\in\g$, $f\in\Md$ and $v\in M$. Then $(\Md,l^*,-l^*-r^*)$ is a
representation of $\g$, which is called the \emph{dual representation}
of $(M,l,r)$. In particular, if $\g$ is free of finite rank, the dual
representation of the adjoint representation $(\g,\ad^L,\ad^R)$ is
$(\g_c^*,\ad^{L*},-\ad^{L*}-\ad^{R*})$.
\end{exam}

Recall that $\Md\otimes M=\Md\otimes_{\C}M$
is a $\C[\partial]$-module via $\partial\otimes1+1\otimes\partial$, and
that $\CEnd(M)$ is a $\C[\partial]$-module via
$(\partial A)_\lambda=-\lambda A_\lambda$. For $f\in\Md$, $u\in M$, write
$\langle f,v\rangle_\nu=\sum_kc_k(v)\nu^k$ (a finite sum, $c_k(v)\in\C$) and set
\[
\Upsilon(f\otimes u)_\lambda(v)=\langle f,v\rangle_{\lambda+\partial}\,u
=\sum_kc_k(v)(\lambda+\partial)^ku\in M[\lambda],
\]
where $\partial$ acts on $u$. We check that this defines a
$\C[\partial]$-linear map $\Upsilon:\Md\otimes M\to\CEnd(M)$.

\emph{(a) Well defined on the tensor product.} The expression
$\langle f,v\rangle_{\lambda+\partial}u$ is $\C$-linear in $f$, in $u$ and in $v$.
Since the tensor product is taken over $\C$, the $\C$-bilinear map
$(f,u)\mapsto\langle f,\cdot\rangle_{\lambda+\partial}u$ factors uniquely through
$\Md\otimes_{\C}M$; no further relation has to be checked.

\emph{(b) Conformal linearity.} Since $\langle f,\partial v\rangle_\nu=\nu\langle f,v\rangle_\nu$,
\[
\Upsilon(f\otimes u)_\lambda(\partial v)=(\lambda+\partial)\langle f,v\rangle_{\lambda+\partial}u
=(\partial+\lambda)\Upsilon(f\otimes u)_\lambda(v),
\]
so $\Upsilon(f\otimes u)\in\CEnd(M)$.

\emph{(c) $\C[\partial]$-linearity.} Since $\langle\partial f,v\rangle_\nu=-\nu\langle f,v\rangle_\nu$,
\[
\Upsilon(\partial f\otimes u+f\otimes\partial u)_\lambda(v)
=-(\lambda+\partial)\langle f,v\rangle_{\lambda+\partial}u+\partial\langle f,v\rangle_{\lambda+\partial}u
=-\lambda\,\Upsilon(f\otimes u)_\lambda(v),
\]
that is, $\Upsilon(\partial(f\otimes u))=\partial\,\Upsilon(f\otimes u)$.

\emph{(d) Bijectivity when $M$ is free of finite rank.} Let
$e_1,\ldots,e_N$ be a basis of $M$. Since $M$ is free, an element
$f\in\Md$ is determined by the polynomials $f_\lambda(e_i)\in\C[\lambda]$,
which can be chosen arbitrarily. Hence $\Md$ is free with dual basis
$e_1^*,\ldots,e_N^*$, $\langle e_i^*,e_j\rangle_\lambda=\delta_{ij}$, and
$\langle p(\partial)e_i^*,e_j\rangle_\lambda=\delta_{ij}p(-\lambda)$. Consequently,
\[
\Md\otimes M=\bigoplus_{i,j=1}^N\C[\partial]e_i^*\otimes\C[\partial]e_j
\cong\bigoplus_{i,j=1}^N\C[s,t],
\qquad p(\partial)e_i^*\otimes q(\partial)e_j\longmapsto p(s)q(t).
\]
Similarly, by conformal linearity, $A\in\CEnd(M)$ is determined by
$A_\lambda(e_k)=\sum_jP_{kj}(\lambda,\partial)e_j$, where the polynomials
$P_{kj}\in\C[\lambda,\partial]$ are arbitrary, so
$\CEnd(M)\cong\bigoplus_{k,j}\C[\lambda,\partial]$. In these coordinates,
\[
\Upsilon(p(\partial)e_i^*\otimes q(\partial)e_j)_\lambda(e_k)=\delta_{ik}\,p(-\lambda-\partial)q(\partial)e_j,
\]
so $\Upsilon$ maps the $(i,j)$-summand to the $(i,j)$-summand by the
algebra homomorphism $\C[s,t]\to\C[\lambda,\partial]$, $s\mapsto-\lambda-\partial$,
$t\mapsto\partial$. This homomorphism is an isomorphism, with inverse
$\lambda\mapsto-s-t$, $\partial\mapsto t$. Hence $\Upsilon$ is injective and
surjective.

Thus, for $M$ free of finite rank, $\Upsilon$ is an isomorphism of
$\C[\partial]$-modules, and we identify $f\otimes u\in\Md\otimes M$ with the
conformal linear map $A_\lambda(v)=\langle f,v\rangle_{\lambda+\partial}u$. Whenever an
element of $\Md\otimes M$ is regarded as an element of $\CEnd(M)$, or
conversely, this identification is meant.

We now show that, under this identification, the right
$\lambda$-action becomes a $1$-cocycle.

\begin{proposition}\label{prop:dual}
Let $(M,l,r)$ be a representation of $\g$, where $M$ is a free
$\C[\partial]$-module of finite rank. Then:
\begin{enumerate}
\item[(i)] for every $y\in\g$, the map $r(y)$ given by
$r(y)_\mu(v)=v_{-\mu-\partial}y$ belongs to $\CEnd(M)$, and
$r:\g\to\CEnd(M)\cong\Md\otimes M$ is $\C[\partial]$-linear;
\item[(ii)] for all $x\in\g$ and $A\in\Md\otimes M\cong\CEnd(M)$,
\[
(\LL)(x)_\lambda A=[l(x)_\lambda A];
\]
\item[(iii)] $\widetilde r=\Upsilon^{-1}\circ r:\g\to\Md\otimes M$ is a
$1$-cocycle of $\g$ with coefficients in the representation
$(\Md\otimes M,\LL,0)$; equivalently,
$r([x_\lambda y])=[l(x)_\lambda r(y)]$ for all $x,y\in\g$, coefficientwise
in $\lambda$.
\end{enumerate}
\end{proposition}

\begin{proof}
(i) Fix $y\in\g$ and $v\in M$, and write
$F(\nu)=v_\nu y\in M[\nu]$. By the right sesquilinearity
identities \eqref{eq:rep-sesq-right},
\[
(\partial v)_\nu y=-\nu F(\nu),
\qquad
v_\nu(\partial y)=(\partial+\nu)F(\nu).
\]
Substituting $\nu=-\mu-\partial$ gives
\[
r(y)_\mu(\partial v)
=(\partial+\mu)r(y)_\mu(v),
\]
and
\[
r(\partial y)_\mu(v)
=-\mu r(y)_\mu(v)
=(\partial r(y))_\mu(v).
\]
Thus $r(y)\in\CEnd(M)$ for every $y\in\g$, and
$r:\g\to\CEnd(M)$ is $\C[\partial]$-linear.

(ii) Write $A=f\otimes u$. Recall that in $\langle f,v\rangle_{\mu+\partial}\,u$ the
operator $\partial$ acts on $u$. By the identification $\Upsilon$,
\[
\begin{aligned}
\big((\LL)(x)_\lambda A\big)_\mu(v)
&=\langle l^*(x)_\lambda f,v\rangle_{\mu+\partial}\,u+\langle f,v\rangle_{\mu+\partial}\,x_\lambda u\\
&=-\langle f,x_\lambda v\rangle_{\mu-\lambda+\partial}\,u+\langle f,v\rangle_{\mu+\partial}\,x_\lambda u .
\end{aligned}
\]
On the other hand, by \eqref{eq:cend-bracket},
\[
[l(x)_\lambda A]_\mu(v)=x_\lambda\big(A_{\mu-\lambda}(v)\big)-A_{\mu-\lambda}(x_\lambda v)
=x_\lambda\big(\langle f,v\rangle_{\mu-\lambda+\partial}u\big)-\langle f,x_\lambda v\rangle_{\mu-\lambda+\partial}\,u .
\]
Since $x_\lambda(\partial w)=(\partial+\lambda)x_\lambda w$, each factor
$\mu-\lambda+\partial$ becomes $\mu+\partial$ after applying $x_\lambda$,
so 
$$x_\lambda(\langle f,v\rangle_{\mu-\lambda+\partial}u)=\langle f,v\rangle_{\mu+\partial}\,x_\lambda u.$$
Hence $$(\LL)(x)_\lambda A=[l(x)_\lambda A].$$

(iii) By (i), $r\in C^1_{\cLeib}(\g,\Md\otimes M)$. Since the right
action of $(\Md\otimes M,\LL,0)$ is zero, \eqref{eq:delta} and (ii) give
\[
(\delta r)_\lambda(x,y)=(\LL)(x)_\lambda r(y)-r([x_\lambda y])
=[l(x)_\lambda r(y)]-r([x_\lambda y]).
\]
For $x,y\in\g$ and $v\in M$, we have
\[
[l(x)_\lambda r(y)]_\mu(v)
=
x_\lambda\big(r(y)_{\mu-\lambda}(v)\big)
-
r(y)_{\mu-\lambda}(x_\lambda v).
\]
We compute the two terms separately, using the convention of
\eqref{eq:endo-convention} for the substitution.

\emph{First term.}
By definition,
 $r(y)_{\mu-\lambda}(v) = v_{-\mu+\lambda-\partial}y$. Writing $v_\nu y=\sum_{k}\nu^k g_k$ (a finite sum), we have
$$r(y)_{\mu-\lambda}(v)=\sum_k(\lambda-\mu-\partial)^k g_k.$$ Applying the
left action of $x$ and using
$x_\lambda(\partial w)=(\partial+\lambda)x_\lambda w$, each factor
$\lambda-\mu-\partial$ becomes $\lambda-\mu-(\partial+\lambda)=-\mu-\partial$,
so that
\[
x_\lambda\big(r(y)_{\mu-\lambda}(v)\big)
=
\left.
x_\lambda(v_\nu y)
\right|_{\nu=-\mu-\partial}.
\]
In this expression, the substitution is performed after applying
$x_\lambda$; thus $\partial$ acts on the resulting coefficients.

\emph{Second term.}
By the definition of $r(y)$,
 $r(y)_{\mu-\lambda}(x_\lambda v) = (x_\lambda v)_{\lambda-\mu-\partial}y$.

\emph{Comparison.}
Identity \eqref{eq:rep2}, with $\nu$ in place of $\mu$, gives
\[
x_\lambda(v_\nu y)
-
(x_\lambda v)_{\lambda+\nu}y
=
v_\nu[x_\lambda y].
\]
Since the right $\lambda$-action is polynomial, the substitution
$\nu=-\mu-\partial$ is defined coefficientwise. It follows that
\[
\left.
x_\lambda(v_\nu y)
\right|_{\nu=-\mu-\partial}
-
(x_\lambda v)_{\lambda-\mu-\partial}y
=
v_{-\mu-\partial}[x_\lambda y].
\]
The right-hand side is precisely
$r([x_\lambda y])_\mu(v)$. 
Combining the preceding identities, we obtain
\[
[l(x)_\lambda r(y)]_\mu(v)
=
r([x_\lambda y])_\mu(v)
\]
for all $x,y\in\g$ and $v\in M$. Hence
\[
[l(x)_\lambda r(y)]
=
r([x_\lambda y]),
\]
where the right-hand side is interpreted coefficientwise in $\lambda$. Therefore $(\delta r)_\lambda(x,y)=0$, and $r$ is a $1$-cocycle of $\g$
with coefficients in $(\Md\otimes M,\LL,0)$.
\end{proof}

\section{Omni-representations of Leibniz conformal algebras}
\label{sec:omni-represent}

In this section, we introduce omni-representations of Leibniz
conformal algebras. Since a representation of a Leibniz conformal
algebra involves both left and right $\lambda$-actions, it cannot in
general be fully encoded by a homomorphism into $\CEnd(M)$ alone.
We therefore use the conformal omni-Lie algebra $\olc(M)$
\cite{TZhang} as an ambient Leibniz conformal algebra.

Let $M$ be
a $\C[\partial]$-module. Then the triple $(M,\Id,0)$ is a
representation of the (formal) Lie conformal algebra $\CEnd(M)$. Indeed, \eqref{eq:rep1} is precisely
\[
A_\lambda(B_\mu v)
=
[A_\lambda B]_{\lambda+\mu}(v)
+
B_\mu(A_\lambda v),
\]
while the identities involving the right action are trivially satisfied.
We may therefore form the corresponding semidirect product.

\begin{definition}[\cite{TZhang}]\label{def:olc}
The \emph{conformal omni-Lie algebra} of $M$ is the semidirect product
$\olc(M)=\CEnd(M)\ltimes_{(\Id,0)}M$, that is, $\CEnd(M)\oplus M$
equipped with the $\lambda$-bracket
\begin{equation}\label{eq:olc}
\{(A,v)_\lambda(B,w)\}
=
\big([A_\lambda B],A_\lambda(w)\big).
\end{equation}
\end{definition}

As a semidirect product, $\olc(M)$ is a formal Leibniz conformal
algebra, and a Leibniz conformal algebra when $M$ is finitely generated
(see the discussion after \eqref{eq:cend-bracket}). In general,
$\olc(M)$ is not a Lie conformal algebra, since
$\{(0,v)_\lambda(A,0)\}=0$,
whereas
$\{(A,0)_\lambda(0,v)\}=(0,A_\lambda(v))$
is nonzero in general.

\emph{Convention on brackets.} Throughout the paper, braces
$\{\cdot_\lambda\cdot\}$ denote only the $\lambda$-bracket of $\olc(M)$.
Square brackets $[\cdot_\lambda\cdot]$ denote the $\lambda$-bracket of $\g$
(or of another Leibniz conformal algebra $\g'$, then written
$[\cdot_\lambda\cdot]'$) and the conformal commutator
\eqref{eq:cend-bracket} of $\CEnd(M)$; the meaning is always determined
by the arguments. The bracket induced on $M$ by a map $T$ is written
$[\cdot_\lambda\cdot]_T$, and juxtaposition $x_\lambda v$, $v_\lambda x$
denotes the $\lambda$-actions of a representation.

Subalgebras of $\olc(M)$ can be obtained from graphs. For a
$\C[\partial]$-linear map $T:M\to\CEnd(M)$, consider its graph
$$\mathcal{G}_T=\{(T(v),v)\mid v\in M\},$$ and define
$\Gamma_T:M\to\mathcal{G}_T,v\mapsto(T(v),v).$ Let
$\pr:\CEnd(M)\oplus M\to M$ be the second projection.

The following proposition characterizes the maps $T$ whose graph is a
subalgebra of $\olc(M)$, and shows that such a graph induces a Leibniz
conformal algebra structure on $M$.

\begin{proposition}\label{prop:graph}
The graph $\mathcal{G}_T$ is a subalgebra of
$\olc(M)$ if and only if
\begin{equation}\label{eq:embedding-tensor}
[T(v)_\lambda T(w)]
=
T[[\lambda]]\big(T(v)_\lambda w\big),
\qquad\forall~ v,w\in M.
\end{equation}
If this condition holds, then
$[v_\lambda w]_T=T(v)_\lambda w$
defines a Leibniz conformal algebra structure $M_T$ on $M$, and
$\Gamma_T:M_T\to\mathcal{G}_T$ is an isomorphism of Leibniz conformal
algebras.
\end{proposition}

\begin{proof}
By the definition \eqref{eq:olc} of the conformal omni-Lie bracket,
\[
\{\Gamma_T(v)_\lambda\Gamma_T(w)\}
=
\big([T(v)_\lambda T(w)],T(v)_\lambda w\big).
\]
Here $T[[\lambda]]$ denotes the coefficientwise extension of $T$
introduced after \eqref{eq:cend-bracket}; since $T(v)_\lambda w$ is a
polynomial in $\lambda$ with coefficients in $M$, so is
$T[[\lambda]](T(v)_\lambda w)$. Hence this expression belongs to
$\mathcal{G}_T[[\lambda]]$ if and only if \eqref{eq:embedding-tensor} holds for $v,w$.
Consequently, $\mathcal{G}_T$ is closed under the $\lambda$-bracket of
$\olc(M)$ if and only if \eqref{eq:embedding-tensor} holds.

Under this condition,
$$\{\Gamma_T(v)_\lambda\Gamma_T(w)\}=\Gamma_T([v_\lambda w]_T)$$ with
$[v_\lambda w]_T=T(v)_\lambda w$.
We check directly that $M_T$ is a Leibniz conformal algebra. The
bracket $[v_\lambda w]_T$ is polynomial in $\lambda$, since
$T(v)\in\CEnd(M)$. Since $T$ is $\C[\partial]$-linear and $T(v)$ is
conformal linear,
\[
[(\partial v)_\lambda w]_T=(\partial T(v))_\lambda w=-\lambda[v_\lambda w]_T,
\qquad
[v_\lambda\partial w]_T=(\partial+\lambda)[v_\lambda w]_T .
\]
Moreover, for $u,v,w\in M$, by \eqref{eq:cend-bracket} and
\eqref{eq:embedding-tensor},
\begin{align*}
[u_\lambda[v_\mu w]_T]_T-[v_\mu[u_\lambda w]_T]_T
&=[T(u)_\lambda T(v)]_{\lambda+\mu}w\\
&=T\big(T(u)_\lambda v\big)_{\lambda+\mu}w
=[[u_\lambda v]_T{}_{\lambda+\mu}w]_T ,
\end{align*}
which is the Leibniz conformal identity \eqref{eq:leibniz} for
$[\cdot_\lambda\cdot]_T$. Hence $M_T$ is a Leibniz conformal algebra.
Finally, $\Gamma_T$ is a $\C[\partial]$-linear bijection from $M$ onto
$\mathcal{G}_T$, with inverse $\pr|_{\mathcal{G}_T}$, and it preserves the
$\lambda$-brackets by the identity displayed above. Hence
$\Gamma_T:M_T\to\mathcal{G}_T$ is an isomorphism.
\end{proof}

Note that $M_T$ and $M$ have the same underlying
$\C[\partial]$-module; the notation $M_T$ indicates that $M$ is endowed
with the $\lambda$-bracket $[\cdot_\lambda\cdot]_T$.

\begin{remark}\label{rmk:embedding}
Recall that, given a representation $(M,l,0)$ of $\g$, a
$\C[\partial]$-linear map $T:M\to\g$ is called an \emph{embedding tensor}
if $[T(u)_\lambda T(v)]=T\big(l(T(u))_\lambda v\big)$ for all $u,v\in M$
\cite{KS,STZ}. Thus \eqref{eq:embedding-tensor} is precisely the
embedding tensor condition for the (formal) Lie conformal algebra
$\CEnd(M)$ and its representation $(M,\Id,0)$. We therefore call a map
$T$ satisfying \eqref{eq:embedding-tensor} a \emph{conformal embedding
tensor} on $M$. Since $T([v_\lambda w]_T)=[T(v)_\lambda T(w)]$, the map
$T:M_T\to\CEnd(M)$ is a homomorphism of (formal) Leibniz conformal
algebras.
\end{remark}

We now introduce the main notion of this paper.

\begin{definition}\label{def:omni-representation}
An \emph{omni-representation} of a Leibniz conformal algebra $\g$ on a
$\C[\partial]$-module $M$ is a $\C[\partial]$-linear map
$\rho:\g\to\olc(M)$ such that
\[
\rho([x_\lambda y])=\{\rho(x)_\lambda\rho(y)\},
\qquad\forall~ x,y\in\g.
\]
\end{definition}

Thus an omni-representation is a homomorphism of (formal) Leibniz
conformal algebras $\g\to\olc(M)$, the identity being understood in
$\olc(M)[[\lambda]]$. In particular, $\im(\rho)$ is closed under the
$\lambda$-bracket of $\olc(M)$, and it is a Leibniz conformal
algebra.

Every $\C[\partial]$-linear map $\rho:\g\to\olc(M)$ can be uniquely
written as
$\rho=(\varphi,\theta)$,
where $\varphi:\g\to\CEnd(M)$ and $\theta:\g\to M$. The homomorphism
condition then splits into two parts.

The next result characterizes omni-representations in terms of a
representation with zero right $\lambda$-action together with a
$1$-cocycle.

\begin{proposition}\label{prop:omni-characterization}
A $\C[\partial]$-linear map
$\rho=(\varphi,\theta):\g\to\CEnd(M)\oplus M$ is an omni-representation
if and only if $(M,\varphi,0)$ is a representation of $\g$ and
$\theta$ is a $1$-cocycle of $\g$ with coefficients in $(M,\varphi,0)$.
\end{proposition}

\begin{proof}
By the definition \eqref{eq:olc} of the bracket of $\olc(M)$, we have
\begin{align*}
\rho([x_\lambda y])=&\big(\varphi([x_\lambda y]),\theta([x_\lambda y])\big),\\
\{\rho(x)_\lambda\rho(y)\}=&\big([\varphi(x)_\lambda\varphi(y)],\varphi(x)_\lambda\theta(y)\big).
\end{align*}
Hence $\rho$ is a homomorphism if and only if
\begin{align}
\varphi([x_\lambda y])
&=
[\varphi(x)_\lambda\varphi(y)],
\label{eq:omni-1}
\\
\theta([x_\lambda y])
&=
\varphi(x)_\lambda\theta(y).
\label{eq:omni-2}
\end{align}
Since $\varphi$ is $\C[\partial]$-linear with values in
$\CEnd(M)$, the left action $x_\lambda v=\varphi(x)_\lambda v$ satisfies
\eqref{eq:rep-sesq-left}, and the first identity, evaluated at $v\in M$
by means of \eqref{eq:cend-bracket}, is exactly the representation identity
\eqref{eq:rep1} for $(M,\varphi,0)$, while the remaining representation
identities involving the right action are trivially satisfied. For
$n=1$ and zero right action, \eqref{eq:delta} gives
\[
(\delta\theta)_\lambda(x,y)
=
\varphi(x)_\lambda\theta(y)-\theta([x_\lambda y]).
\]
Thus \eqref{eq:omni-2} is equivalent to $\delta\theta=0$.
\end{proof}

We now show that every representation of $\g$ on a free
$\C[\partial]$-module of finite rank naturally gives rise to an
omni-representation. The construction combines
Proposition~\ref{prop:omni-characterization} with the fact, established
in Proposition~\ref{prop:dual}, that the right $\lambda$-action is a
$1$-cocycle with values in $\Md\otimes M\cong\CEnd(M)$.

\begin{thm}\label{thm:induced-omni-rep}
Let $(M,l,r)$ be a representation of $\g$, where $M$ is a free
$\C[\partial]$-module of finite rank. Then
\[
\rho=(\LL,\,\widetilde r):\g\longrightarrow\olc(\Md\otimes M)
=\CEnd(\Md\otimes M)\oplus(\Md\otimes M),
\qquad \widetilde r=\Upsilon^{-1}\circ r,
\]
is an omni-representation of $\g$ on $\Md\otimes M$. Here $r$ takes
values in $\CEnd(M)$, while $\widetilde r$ is its expression in
$\Md\otimes M$. Here
$\olc(\Md\otimes M)$ is a formal Leibniz conformal algebra, since
$\Md\otimes M$ is not finitely generated, and the homomorphism identity
$\rho([x_\lambda y])=\{\rho(x)_\lambda\rho(y)\}$ is a priori an identity in
$\olc(\Md\otimes M)[[\lambda]]$; both sides are in fact polynomial in
$\lambda$. In particular, the image $\im(\rho)$ is a genuine Leibniz
conformal algebra and, by Lemma~\ref{lem:image-rep}, a genuine
representation of $\g$.
\end{thm}

\begin{proof}
Hypotheses (H1) and (H2) hold, since $M$ is free of finite rank. By
(H1), $(\Md\otimes M,\LL,0)$ is a representation of $\g$; in particular
$(\LL)(x)\in\CEnd(\Md\otimes M)$ for every $x\in\g$. By (H2), $\Upsilon$ is
an isomorphism, and by Proposition~\ref{prop:dual},
$\widetilde r=\Upsilon^{-1}\circ r:\g\to\Md\otimes M$ is a $1$-cocycle on $\g$ with
coefficients in $(\Md\otimes M,\LL,0)$. This cocycle identity is
Proposition~\ref{prop:dual}(iii), whose proof uses only the
coefficientwise substitution convention \eqref{eq:endo-convention},
applied to polynomials. By
Proposition~\ref{prop:omni-characterization}, $\rho=(\LL,\,\widetilde r)$ is a
homomorphism from $\g$ to $\olc(\Md\otimes M)$, that is, an
omni-representation of $\g$ on $\Md\otimes M$.

Note that, by part~(d) of the construction of $\Upsilon$,
$\Md\otimes M\cong\bigoplus_{i,j}\C[s,t]$ with $\partial$ acting as
multiplication by $s+t$; this module is free of infinite rank over
$\C[\partial]$. Hence $\CEnd(\Md\otimes M)$ is only a formal Lie
conformal algebra, and $\olc(\Md\otimes M)$ a formal Leibniz conformal
algebra. Finally, the two
components of the homomorphism identity are
$[(\LL)(x)_\lambda(\LL)(y)]=(\LL)([x_\lambda y])$ and
$(\LL)(x)_\lambda\widetilde r(y)=\widetilde r([x_\lambda y])$, whose right-hand sides are
polynomial in $\lambda$ because the $\lambda$-bracket of $\g$ is.
\end{proof}

The assumption that $M$ is free of finite rank in
Theorem~\ref{thm:induced-omni-rep} is used to identify $\CEnd(M)$ with
$\Md\otimes M$, and this identification fails without it: if
$p(\partial)v=0$ with $p\neq0$, then $p(\lambda)\langle f,v\rangle_\lambda=0$,
so $\langle f,v\rangle_\lambda=0$ for all $f\in\Md$; thus $\Md$ does not see
the torsion of $M$, and $\Upsilon$ is not an isomorphism in general.
The omni-representation $\rho$
contains the right action itself, through $\widetilde r$, whereas the left action appears
only through $\LL$.
\begin{exam}\label{ex:current-omni}
Let $\g=\Cur\lie=\C[\partial]e_1\oplus\C[\partial]e_2$ be the Leibniz
conformal algebra of Example~\ref{ex:current}, and let $M=\g$ be its
adjoint representation. Recall that the only nonzero bracket between the
generators is $[e_{2\,\lambda}e_2]=e_1$. Let $e_1^*,e_2^*\in\g_c^*$ be the
dual basis, so that $\langle e_i^*,e_j\rangle_\lambda=\delta_{ij}$.

For $v=q(\partial)e_2+s(\partial)e_1\in\g$, sesquilinearity gives
\[
e_{2\,\lambda}v=q(\partial+\lambda)e_1,
\qquad
v_\lambda e_2=q(-\lambda)e_1,
\]
and all the other actions of the generators vanish. Substituting
$\lambda=-\mu-\partial$ as in \eqref{eq:endo-convention}, we find
$r(e_2)_\mu(v)=q(\partial+\mu)e_1=l(e_2)_\mu(v)$, while $r(e_1)=l(e_1)=0$.

Since $\langle e_2^*,v\rangle_\nu=q(\nu)$, we have
$\Upsilon(e_2^*\otimes e_1)_\mu(v)=q(\partial+\mu)e_1=r(e_2)_\mu(v)$, that is,
$\widetilde r(e_2)=\Upsilon^{-1}(r(e_2))=e_2^*\otimes e_1$. By
Theorem~\ref{thm:induced-omni-rep}, the omni-representation of $\g$ on
$\g_c^*\otimes\g$ is therefore
\[
\rho(e_1)=0,\qquad
\rho(e_2)=\big((\LL)(e_2),\,e_2^*\otimes e_1\big),
\]
and its cocycle component is nonzero. One can also check the cocycle
condition directly: since $e_{2\,\lambda}v$ is a multiple of $e_1$,
\[
\langle l^*(e_2)_\lambda e_2^*,v\rangle_\mu
=-\langle e_2^*,q(\partial+\lambda)e_1\rangle_{\mu-\lambda}=0,
\]
and $e_{2\,\lambda}e_1=0$, so that
$(\delta\widetilde r)_\lambda(e_2,e_2)
=(\LL)(e_2)_\lambda(e_2^*\otimes e_1)-\widetilde r(e_1)=0$.
Finally, $\im(\rho)=\C[\partial]\rho(e_2)$ is free of rank one.
\end{exam}

An omni-representation $\rho$ induces a representation of $\g$
on its image $\im(\rho)$. This representation is the coefficient
module of the omni-cohomology in Section~\ref{sec:omni-cohomo}. It is
a special case of the following general construction.

\begin{lemma}\label{lem:pullback}
Let $\theta:\g\to\g'$ be a homomorphism of Leibniz conformal
algebras. Then $\g'$ is a representation of $\g$ with actions
\begin{equation}\label{eq:pullback-actions}
x_\lambda a=[\theta(x)_\lambda a],
\qquad
a_\lambda x=[a_\lambda\theta(x)],
\end{equation}
and $\theta(\g)$ is a subrepresentation of $\g'$.
\end{lemma}

\begin{proof}
The sesquilinearity conditions follow from those in $\g'$ and
from the $\C[\partial]$-linearity of $\theta$. Applying the Leibniz
conformal identity in $\g'$ to the triples
$(\theta(x),\theta(y),a)$, $(\theta(x),a,\theta(y))$ and
$(a,\theta(x),\theta(y))$ gives respectively \eqref{eq:rep1},
\eqref{eq:rep2} and \eqref{eq:rep3}. Moreover,
$$x_\lambda\theta(y)=\theta([x_\lambda y])\quad \text{and}
\quad\theta(y)_\lambda x=\theta([y_\lambda x]),$$ so $\theta(\g)$ is stable under both actions.
\end{proof}

Let $\rho:\g\to\olc(M)$ be an omni-representation. For
$U\in\im(\rho)$, define
\begin{equation}\label{eq:im-actions}
x_\lambda U=\{\rho(x)_\lambda U\},
\qquad
U_\lambda x=\{U_\lambda\rho(x)\}.
\end{equation}

\begin{lemma}\label{lem:image-rep}
With the above notation, $\im(\rho)$ is a $\C[\partial]$-submodule of
$\olc(M)$ and a representation of $\g$. More explicitly, if
$U=\rho(u)$, then $x_\lambda U=\rho([x_\lambda u])$ and $U_\lambda x=\rho([u_\lambda x])$.
In particular, $\rho:\g\to\im(\rho)$ is a surjective homomorphism of
representations from the adjoint representation of $\g$.
\end{lemma}

\begin{proof}
Since $\rho$ is $\C[\partial]$-linear, $\im(\rho)$ is a
$\C[\partial]$-submodule. Since $\rho$ preserves the $\lambda$-brackets, we have
$\{\rho(x)_\lambda\rho(u)\}=\rho([x_\lambda u])$ and
$\{\rho(u)_\lambda\rho(x)\}=\rho([u_\lambda x])$.
These expressions depend only on $U=\rho(u)$ and lie in
$\im(\rho)[\lambda]$, since the $\lambda$-bracket of $\g$ is polynomial. The sesquilinearity conditions and the
identities \eqref{eq:rep1}--\eqref{eq:rep3} then follow by applying
$\rho$ to \eqref{eq:sesq} and to the Leibniz conformal identity
\eqref{eq:leibniz} in $\g$. The last assertion is a restatement of the
two displayed identities.
\end{proof}

We now consider omni-representations whose image lies in the
graph of a conformal embedding tensor.

Let $T:M\to\CEnd(M)$ be a conformal embedding tensor, and let $M_T$ be
the Leibniz conformal algebra of Proposition~\ref{prop:graph}. An
omni-representation $\rho$ is said to \emph{take values in the graph}
$\mathcal{G}_T$ if $\im(\rho)\subseteq\mathcal{G}_T$.

\begin{proposition}\label{prop:graph-omni}
Let $T:M\to\CEnd(M)$ be a conformal embedding tensor. Then the
assignment $\theta\mapsto\rho_\theta=\Gamma_T\circ\theta=(T\circ\theta,\theta)$
is a bijection between homomorphisms
$\theta:\g\to M_T$ and omni-representations of $\g$ on $M$ with values
in $\mathcal{G}_T$. Moreover,
$$\im(\rho_\theta)=\Gamma_T(\im(\theta)).$$
\end{proposition}

\begin{proof}
By Proposition~\ref{prop:graph}, $\mathcal{G}_T$ is a subalgebra of
$\olc(M)$ and $\Gamma_T:M_T\to\mathcal{G}_T$ is an isomorphism of
Leibniz conformal algebras, with inverse
$\pr|_{\mathcal{G}_T}$. Hence a $\C[\partial]$-linear map
$\rho:\g\to\olc(M)$ with values in $\mathcal{G}_T$ is a homomorphism
if and only if $\theta=\pr\circ\rho$ is a homomorphism
$\g\to M_T$, and then $\rho=\Gamma_T\circ\theta$. Since $\Gamma_T$
is a bijection from $M$ onto $\mathcal{G}_T$, we have
$\im(\rho_\theta)=\Gamma_T(\im(\theta))$.
\end{proof}

\begin{remark}\label{rmk:graph-vs-characterization}
For $\rho_\theta=(T\circ\theta,\theta)$, condition \eqref{eq:omni-2} is the
homomorphism identity $\theta([x_\lambda y])=T(\theta(x))_\lambda\theta(y)$, and
\eqref{eq:omni-1} follows from it by applying $T$ and using
\eqref{eq:embedding-tensor}.
\end{remark}

For omni-representations with values in $\mathcal{G}_T$, the
representation of $\g$ on $\im(\rho)$ given by
Lemma~\ref{lem:image-rep} can be transported to $M$ itself. Let
$\rho=\rho_\theta=(T\circ\theta,\theta)$. Define actions on $M$ by
\begin{align}
x_\lambda v
&=
\pr\{\rho(x)_\lambda\Gamma_T(v)\}
=
T(\theta(x))_\lambda v
=
[\theta(x)_\lambda v]_T,
\label{eq:graph-left}
\\
v_\lambda x
&=
\pr\{\Gamma_T(v)_\lambda\rho(x)\}
=
T(v)_\lambda\theta(x)
=
[v_\lambda\theta(x)]_T.
\label{eq:graph-right}
\end{align}

\begin{lemma}\label{lem:graph-rep}
With these actions, $M$ is a representation of $\g$, and $\im(\theta)$
is a subrepresentation. Moreover,
\begin{equation}\label{eq:graph-intertwine}
\{\rho(x)_\lambda\Gamma_T(v)\}
=
\Gamma_T(x_\lambda v),
\qquad
\{\Gamma_T(v)_\lambda\rho(x)\}
=
\Gamma_T(v_\lambda x).
\end{equation}
The restriction
$\Gamma_T:\im(\theta)\to\im(\rho)$
is an isomorphism of representations.
\end{lemma}

\begin{proof}
Since $\theta:\g\to M_T$ is a homomorphism, the actions
\eqref{eq:graph-left} and \eqref{eq:graph-right} are precisely the
pullback actions of Lemma~\ref{lem:pullback}. Thus $M$ is a
representation and $\im(\theta)$ is a subrepresentation. Furthermore,
$$\{\rho(x)_\lambda\Gamma_T(v)\}=\Gamma_T([\theta(x)_\lambda v]_T)=\Gamma_T(x_\lambda v)$$
and
$$\{\Gamma_T(v)_\lambda\rho(x)\}=\Gamma_T([v_\lambda\theta(x)]_T)=\Gamma_T(v_\lambda x).$$
Restricting to $v\in\im(\theta)$ proves the final assertion.
\end{proof}

\begin{remark}\label{rmk:M-vs-image}
By Lemma~\ref{lem:graph-rep}, the representation $\im(\rho)$ is
isomorphic to $\im(\theta)$, and not to $M$ unless $\theta$ is surjective.
For instance, for $M=\C$, $T=0$ and $\theta=0$, we have $\im(\rho_0)=0$,
whereas $M$ is the trivial representation. For a nonzero example, let
$\g$ be finitely generated, let $\C_0$ be $\C$ with $\partial=0$ and zero
bracket, and take $M=\g\oplus\C_0$, $T=\ad^L_M$ and $\theta(x)=(x,0)$,
which is a homomorphism $\g\to M_T=M$ by Lemma~\ref{lem:adL-embedding}
below. Then $M$ is the direct sum of the adjoint and the trivial
representations, whereas $\im(\rho_\theta)\cong\im(\theta)\cong\g$ is the
adjoint representation.
\end{remark}

We now consider two particular cases: $T=0$ on $M=\C$, and
$T=\ad^L_{\g'}$ on a Leibniz conformal algebra $\g'$, which for
$\g'=\g$ gives $T=\ad^L$. We begin with the trivial case.

Let $M=\C$, with $\partial$ acting trivially. If $A\in\CEnd(\C)$, then
\[
0=A_\lambda(\partial 1)
=(\partial+\lambda)A_\lambda(1)
=\lambda A_\lambda(1),
\]
and hence $A_\lambda(1)=0$. Thus
\[
\CEnd(\C)=0,
\qquad
\olc(\C)=\{0\}\oplus\C.
\]
Consequently, $T=0$ is the only $\C[\partial]$-linear map
$\C\to\CEnd(\C)$. Its graph is
$\mathcal{G}_0=\olc(\C)$, and the corresponding Leibniz conformal
algebra $M_T$ is $\C$ with zero $\lambda$-bracket. We denote it by
$\C_0$. An omni-representation of $\g$ on $\C$ is called a
\emph{trivial omni-representation}. Since $\CEnd(\C)=0$, every trivial
omni-representation has the form
\[
\rho_\xi=(0,\xi),
\]
where $\xi:\g\to\C$ is $\C[\partial]$-linear.

\begin{corollary}\label{cor:trivial-omni}
Trivial omni-representations of $\g$ are in one-to-one correspondence
with $\C$-linear maps $\xi:\g\to\C$ such that
\[
\xi(\partial x)=0
\qquad\text{and}\qquad
\xi(x_{(j)}y)=0
\qquad\text{for all }x,y\in\g\text{ and }j\geq0.
\]
The associated representation on $\C$ is the trivial representation.
\end{corollary}

\begin{proof}
By Proposition~\ref{prop:graph-omni}, trivial omni-representations
correspond to homomorphisms $\xi:\g\to\C_0$. Since $\partial$ acts by
zero on $\C$, $\C[\partial]$-linearity is equivalent to
$\xi(\partial x)=0$ for all $x\in\g$. Since the bracket of $\C_0$
is zero, the homomorphism condition is equivalent to
$\xi([x_\lambda y])=0$ for all $x,y\in\g$, that is,
$\xi(x_{(j)}y)=0$ for all $j\geq0$. The last assertion follows from
$T=0$.
\end{proof}

\begin{lemma}\label{lem:adL-embedding}
Let $\mathcal A$ be a Leibniz conformal algebra with $\lambda$-bracket
$[\cdot_\lambda\cdot]$. Then the map
$\ad^L_{\mathcal A}:\mathcal A\to\CEnd(\mathcal A)$,
$\ad^L_{\mathcal A}(u)_\lambda(w)=[u_\lambda w]$, is a conformal
embedding tensor, and
$\mathcal A_{\ad^L_{\mathcal A}}=\mathcal A$. In particular, for
$\mathcal A=\g$, the map $\ad^L=\ad^L_{\g}$ is a conformal embedding
tensor and $\g_{\ad^L}=\g$.
\end{lemma}

\begin{proof}
Write $\ad^L_{\mathcal A}=\ad$ in this proof.
By the conformal sesquilinearity identities, $\ad(u)$ is a conformal
linear map for every $u\in\mathcal A$, and $\ad$ is
$\C[\partial]$-linear. Let $u,w,z\in\mathcal A$. By the definition
of the conformal commutator on $\CEnd(\mathcal A)$,
\[
[\ad(u)_\lambda\ad(w)]_\nu(z)
=
[u_\lambda[w_{\nu-\lambda}z]]
-
[w_{\nu-\lambda}[u_\lambda z]].
\]
Each conformal bracket of fixed elements is polynomial in its formal
variable. Hence both terms above are well-defined elements of
$\mathcal A[\lambda,\nu]$, and the substitution $\mu=\nu-\lambda$
in the Leibniz conformal identity is legitimate. It yields
\[
[u_\lambda[w_{\nu-\lambda}z]]-[w_{\nu-\lambda}[u_\lambda z]]
=[[u_\lambda w]_\nu z]=\ad([u_\lambda w])_\nu(z).
\]
Therefore $[\ad(u)_\lambda\ad(w)]=\ad\big(\ad(u)_\lambda w\big)$, so the
embedding tensor identity \eqref{eq:embedding-tensor} holds for
$T=\ad^L_{\mathcal A}$. Moreover, $[u_\lambda w]_{\ad}=\ad(u)_\lambda w=[u_\lambda w]$,
which proves $\mathcal A_{\ad^L_{\mathcal A}}=\mathcal A$.
\end{proof}

Proposition~\ref{prop:graph-omni} now gives the following.

\begin{corollary}\label{cor:hom-omni}
Let $\theta:\g\to\g'$ be a homomorphism of Leibniz conformal
algebras. Then $\rho_\theta=(\ad^L_{\g'}\circ\theta,\theta):\g\to\olc(\g')$
is an omni-representation. Every omni-representation of $\g$ on
$\g'$ with values in $\mathcal{G}_{\ad^L_{\g'}}$ is of this form. The
associated representation of $\g$ on $\g'$ is
\[
x_\lambda v=[\theta(x)_\lambda v]',
\qquad
v_\lambda x=[v_\lambda\theta(x)]'.
\]
\end{corollary}

\begin{proof}
Apply Proposition~\ref{prop:graph-omni} to the conformal embedding
tensor $\ad^L_{\g'}$ given by Lemma~\ref{lem:adL-embedding} with
$\mathcal A=\g'$. Recall that
$\g'_{\ad^L_{\g'}}=\g'$.
\end{proof}

\begin{proposition}\label{prop:T-is-ad}
Let $M$ be a $\C[\partial]$-module. The assignment $T\mapsto M_T$ is a
bijection from the set of conformal embedding tensors on $M$ onto the
set of Leibniz conformal algebra structures on $M$, whose inverse sends
a structure to its left adjoint operator. In particular,
\[
T=\ad^L_{M_T}
\qquad\text{and}\qquad
\mathcal{G}_T=\mathcal{G}_{\ad^L_{M_T}},
\]
so that Proposition~\ref{prop:graph-omni} for $T$ is
Corollary~\ref{cor:hom-omni} for $\g'=M_T$.
\end{proposition}

\begin{proof}
Let $T:M\to\CEnd(M)$ be a conformal embedding tensor. By the
definition of the $\lambda$-bracket of $M_T$ in
Proposition~\ref{prop:graph},
$$\ad^L_{M_T}(v)_\lambda(w)=[v_\lambda w]_T=T(v)_\lambda w,\quad \text{for all}~ v,w\in M,$$
so that $$T=\ad^L_{M_T}.$$ Conversely, let $[\cdot_\lambda\cdot]'$ be a
Leibniz conformal algebra structure on $M$, and denote by $\mathcal A$
the module $M$ endowed with it. By Lemma~\ref{lem:adL-embedding},
$\ad^L_{\mathcal A}$ is a conformal embedding tensor and
$M_{\ad^L_{\mathcal A}}=\mathcal A$. The two assignments are therefore
inverse to each other. The last assertion follows from
$\mathcal{G}_T=\mathcal{G}_{\ad^L_{M_T}}$.
\end{proof}

\begin{remark}\label{rmk:T-is-ad}
Consequently, Corollary~\ref{cor:hom-omni} and
Proposition~\ref{prop:graph-omni} are equivalent: the omni-representations
with values in a graph are exactly the maps $(\ad^L_{\g'}\circ\theta,\theta)$
induced by homomorphisms $\theta:\g\to\g'$. This includes the trivial
case, since $0=\ad^L_{\C_0}$. Proposition~\ref{prop:T-is-ad} is the Leibniz
conformal counterpart of Weinstein's correspondence between Lie algebra
structures on $V$ and Dirac structures in $\mathfrak{gl}(V)\oplus V$ that are
graphs \cite{Weinstein}.
\end{remark}

Taking $\g'=\g$ and $\theta=\Id_\g$, we obtain the analogue of the
adjoint representation. The \emph{adjoint omni-representation} of a
Leibniz conformal algebra $\g$ on itself is
\[
\Ad=(\ad^L,\Id_{\g}).
\]

\begin{corollary}\label{cor:adjoint-omni}
The map $\Ad$ is an omni-representation,
$\im(\Ad)=\mathcal{G}_{\ad^L}$, and the associated representation of
$\g$ on itself is the adjoint representation.
\end{corollary}

\begin{proof}
Apply Corollary~\ref{cor:hom-omni} with $\g'=\g$ and
$\theta=\Id_{\g}$. Then $$\im(\Ad)=\Gamma_{\ad^L}(\g)=\mathcal{G}_{\ad^L},$$
and the induced actions are
$x_\lambda v=[x_\lambda v]$ and
$v_\lambda x=[v_\lambda x]$.
\end{proof}

\section{Omni-cohomology of Leibniz conformal algebras}
\label{sec:omni-cohomo}

In this section, we introduce the omni-cohomology of a Leibniz conformal
algebra with coefficients in an omni-representation $\rho$. For
omni-representations with values in the graph of a conformal embedding
tensor, we show that the omni-cohomology is isomorphic to the usual
cohomology with coefficients in a suitable subrepresentation, and we
deduce the trivial and adjoint cases.

Let $\rho:\g\to\olc(M)$ be an omni-representation of a Leibniz
conformal algebra $\g$. The image $\im(\rho)\subseteq\olc(M)$ is a
Leibniz conformal subalgebra and, by Lemma~\ref{lem:image-rep}, a
representation of $\g$ with the actions \eqref{eq:im-actions}. We define
\[
C^n_{\omni}(\g;\rho)=
\begin{cases}
\im(\rho)/\partial\im(\rho),&n=0,\\[1mm]
\Homcs\big(\g^{\otimes n},\im(\rho)[[\lambda_1,\ldots,\lambda_{n-1}]]\big),&n\geq1,
\end{cases}
\]
and $\delta_\rho:C^n_{\omni}(\g;\rho)\to C^{n+1}_{\omni}(\g;\rho)$ by
\[
\delta_\rho\big(U+\partial\im(\rho)\big)(x)=-\{U_\lambda\rho(x)\}\big|_{\lambda=0}
\]
in degree zero and, for $n\geq1$, by
\begin{equation}\label{eq:omni-differential}
\begin{aligned}
&(\delta_\rho\Phi)_{\lambda_1,\ldots,\lambda_n}(x_1,\ldots,x_{n+1})\\
&\quad=\sum_{i=1}^{n}(-1)^{i+1}
\big\{\rho(x_i)_{\lambda_i}\,
\Phi_{\lambda_1,\ldots,\widehat{\lambda_i},\ldots,\lambda_n}
(x_1,\ldots,\widehat{x_i},\ldots,x_{n+1})\big\}\\
&\quad+(-1)^{n+1}\big\{\Phi_{\lambda_1,\ldots,\lambda_{n-1}}
(x_1,\ldots,x_n)_{\Lambda_n}\,\rho(x_{n+1})\big\}\\
&\quad+\sum_{1\leq i<j\leq n+1}(-1)^i\,
\Phi_{\lambda_1,\ldots,\widehat{\lambda_i},\ldots,\lambda_{j-1},
\lambda_i+\lambda_j,\ldots,\lambda_n}\\
&\hspace{3em}\big(x_1,\ldots,\widehat{x_i},\ldots,x_{j-1},
[x_{i\,\lambda_i}x_j],x_{j+1},\ldots,x_{n+1}\big).
\end{aligned}
\end{equation}

\begin{lemma}\label{lem:omni-delta}
With the above notation, $\delta_\rho^2=0$.
\end{lemma}

\begin{proof}
By Lemma~\ref{lem:image-rep}, $\im(\rho)$ is a representation of $\g$
with 
$$x_\lambda U=\{\rho(x)_\lambda U\}\quad \text{and} \quad U_\lambda x=\{U_\lambda\rho(x)\},$$
and these actions take values in $\im(\rho)[\lambda]$. Comparing with
\eqref{eq:delta0}--\eqref{eq:delta}, $\delta_\rho$ is the coboundary
operator of $\g$ with coefficients in this representation, so
$\delta_\rho^2=0.$
\end{proof}

\begin{definition}\label{def:omni-cohomology}
It follows from Lemma~\ref{lem:omni-delta} that
$\big(C^\bullet_{\omni}(\g;\rho),\delta_\rho\big)$ is a
\emph{cochain complex}. Let
$Z^n_{\omni}(\g;\rho)=\{\Phi\in C^n_{\omni}(\g;\rho)\mid\delta_\rho\Phi=0\}$
denote the \emph{space of $n$-cocycles} and
$B^n_{\omni}(\g;\rho)=\{\delta_\rho\Psi\mid\Psi\in C^{n-1}_{\omni}(\g;\rho)\}$
the \emph{space of $n$-coboundaries}, with $B^0_{\omni}(\g;\rho)=0$. Since $\delta_\rho^2=0$, we have
$B^n_{\omni}(\g;\rho)\subset Z^n_{\omni}(\g;\rho)$. The corresponding
quotient
\[
H^n_{\omni}(\g;\rho):=\frac{Z^n_{\omni}(\g;\rho)}{B^n_{\omni}(\g;\rho)},
\qquad\text{for } n\geq0,
\]
is called the \emph{$n$-th omni-cohomology group} of $\g$ with
coefficients in the omni-representation $\rho$. By construction,
$H^\bullet_{\omni}(\g;\rho)$ is the cohomology of $\g$ with coefficients
in the representation $\im(\rho)\subseteq\olc(M)$ of
Lemma~\ref{lem:image-rep}, and not in $M$ itself.
\end{definition}

We now consider the case of graphs. Suppose that $\rho$ has its image
contained in a graph $\mathcal{G}_T$, where
$T:M\to\CEnd(M)$ is a conformal embedding tensor. By
Proposition~\ref{prop:graph-omni},
$\rho=\rho_\theta=(T\circ\theta,\theta)$
for a homomorphism $\theta:\g\to M_T$. In this case, every
$n$-cochain $\Phi\in C^n_{\omni}(\g;\rho)$, $n\geq1$, is uniquely
determined by
$f=\pr\circ\Phi\in C^n_{\cLeib}(\g,\im(\theta))$ through
\begin{equation}\label{eq:graph-cochain}
\Phi=\Gamma_T\circ f=(T\circ f,f).
\end{equation}

For $f\in C^n_{\cLeib}(\g,\im(\theta))$, $n\geq1$, and
$x_1,\ldots,x_{n+1}\in\g$, set
\begin{align*}
f^{(i)}
&=
f_{\lambda_1,\ldots,\widehat{\lambda_i},\ldots,\lambda_n}
(x_1,\ldots,\widehat{x_i},\ldots,x_{n+1}),
\quad 1\leq i\leq n,
\\
f^{(0)}
&=
f_{\lambda_1,\ldots,\lambda_{n-1}}(x_1,\ldots,x_n),
\\
f^{(i,j)}
&=
f_{\lambda_1,\ldots,\widehat{\lambda_i},\ldots,
\lambda_{j-1},\lambda_i+\lambda_j,\ldots,\lambda_n}
\\
&\qquad
(x_1,\ldots,\widehat{x_i},\ldots,x_{j-1},
[x_{i\,\lambda_i}x_j],x_{j+1},\ldots,x_{n+1}),
\quad 1\leq i<j\leq n+1.
\end{align*}

With the left and right $\lambda$-actions
\eqref{eq:graph-left} and \eqref{eq:graph-right}, namely
\[
x_\lambda v=T(\theta(x))_\lambda v,
\qquad
v_\lambda x=T(v)_\lambda\theta(x),
\]
the coboundary operator on cochains with coefficients in
$\im(\theta)$ is
\begin{equation}\label{eq:graph-delta}
\begin{aligned}
(\delta f)_{\lambda_1,\ldots,\lambda_n}
(x_1,\ldots,x_{n+1})
={}&
\sum_{i=1}^{n}(-1)^{i+1}
T(\theta(x_i))_{\lambda_i}\big(f^{(i)}\big)
\\
&+
(-1)^{n+1}
T\big(f^{(0)}\big)_{\Lambda_n}
\big(\theta(x_{n+1})\big)
\\
&+
\sum_{1\leq i<j\leq n+1}(-1)^if^{(i,j)}
\end{aligned}
\end{equation}
for $n\geq1$, and
\[
\delta(v+\partial\im(\theta))(x)
=
-T(v)_\lambda\theta(x)\big|_{\lambda=0}
\]
in degree zero.

The following theorem identifies the coefficient representation
$\im(\rho)$ of the omni-cohomology with the subrepresentation
$\im(\theta)$ of $M$. 

\begin{thm}\label{thm:graph-omni-cohomology}
Let $T:M\to\CEnd(M)$ be a conformal embedding tensor, and let $\rho$
be an omni-representation of $\g$ on $M$ with values in
$\mathcal{G}_T$. By Proposition~\ref{prop:graph-omni},
$\rho=\rho_\theta=(T\circ\theta,\theta)$ for a unique homomorphism
$\theta:\g\to M_T$. Then
$f\mapsto\Gamma_T\circ f$ is an isomorphism of cochain complexes
$\big(C^\bullet_{\cLeib}(\g,\im(\theta)),\delta\big)\to\big(C^\bullet_{\omni}(\g;\rho),\delta_\rho\big)$,
where $\delta$ is given by \eqref{eq:graph-delta}. Consequently,
\[
H^\bullet_{\omni}(\g;\rho)\cong H^\bullet(\g,\im(\theta)),
\]
where $H^\bullet(\g,\im(\theta))$ denotes the cohomology of $\g$ with
coefficients in the subrepresentation $\im(\theta)\subseteq M$ of
Lemma~\ref{lem:graph-rep}, whose $\lambda$-actions are given by
\eqref{eq:graph-left} and \eqref{eq:graph-right}.
\end{thm}

\begin{proof}
By Lemma~\ref{lem:graph-rep},
$\Gamma_T:\im(\theta)\to\im(\rho)$
is a $\C[\partial]$-linear bijection with inverse
$\pr|_{\im(\rho)}$. We extend $\Gamma_T$ coefficientwise to formal power series
in $\lambda_1,\ldots,\lambda_{n-1}$ and to conformal cochains. Since $\Gamma_T$ is a
$\C[\partial]$-linear bijection, this extension preserves conformal
sesquilinearity and satisfies
\[
\Gamma_T(\partial\im(\theta))
=
\partial\im(\rho).
\]
Consequently, for every $n\geq0$, $\Gamma_T$ induces a bijection
$C^n_{\cLeib}(\g,\im(\theta))\to C^n_{\omni}(\g;\rho)$. For $n=0$, it is
the isomorphism
\[
\overline{\Gamma}_T:
\frac{\im(\theta)}{\partial\im(\theta)}
\longrightarrow
\frac{\im(\rho)}{\partial\im(\rho)},
\qquad
v+\partial\im(\theta)
\longmapsto
\Gamma_T(v)+\partial\im(\rho).
\]
For $n\geq1$, it is the map $f\mapsto\Gamma_T\circ f$, whose inverse is
$\Phi\mapsto\pr\circ\Phi$ by \eqref{eq:graph-cochain}.

It remains to show that these maps commute with the differentials. The
first identity in \eqref{eq:graph-intertwine} handles the terms
$\{\rho(x_i)_{\lambda_i}\cdot\}$, and the second one handles the term
$\{\cdot_{\Lambda_n}\rho(x_{n+1})\}$. Hence, by \eqref{eq:omni-differential},
\begin{align*}
&(\delta_\rho(\Gamma_T\circ f))_{\lambda_1,\ldots,\lambda_n}(x_1,\ldots,x_{n+1})\\
&\quad=\sum_{i=1}^{n}(-1)^{i+1}\{\rho(x_i)_{\lambda_i}\Gamma_T(f^{(i)})\}
+(-1)^{n+1}\{\Gamma_T(f^{(0)})_{\Lambda_n}\rho(x_{n+1})\}\\
&\qquad+\sum_{1\leq i<j\leq n+1}(-1)^i\,\Gamma_T(f^{(i,j)})\\
&\quad=\Gamma_T\Big(\sum_{i=1}^{n}(-1)^{i+1}x_{i\,\lambda_i}f^{(i)}
+(-1)^{n+1}f^{(0)}{}_{\Lambda_n}x_{n+1}
+\sum_{1\leq i<j\leq n+1}(-1)^if^{(i,j)}\Big)\\
&\quad=\Gamma_T\big((\delta f)_{\lambda_1,\ldots,\lambda_n}(x_1,\ldots,x_{n+1})\big).
\end{align*}
Similarly, for $v\in\im(\theta)$ and $x\in\g$, the second identity in
\eqref{eq:graph-intertwine} gives
\begin{align*}
\delta_\rho\big(\overline{\Gamma}_T(v+\partial\im(\theta))\big)(x)
&=-\{\Gamma_T(v)_\lambda\rho(x)\}\big|_{\lambda=0}
=-\Gamma_T\big(T(v)_\lambda\theta(x)\big)\big|_{\lambda=0}\\
&=\overline{\Gamma}_T\big(\delta(v+\partial\im(\theta))(x)\big).
\end{align*}
Therefore $\Gamma_T$ is an isomorphism of cochain complexes, and the
result follows.

\end{proof}

\begin{remark}\label{rmk:compare-LS}
(1) With the notation of Theorem~\ref{thm:graph-omni-cohomology}, denote
by $H^\bullet(\g,M)$ the cohomology of $\g$ with coefficients in the
whole module $M$, endowed with the actions \eqref{eq:graph-left} and
\eqref{eq:graph-right} of Lemma~\ref{lem:graph-rep}. The theorem shows
that the coefficient module of the omni-cohomology is $\im(\theta)$, and
not $M$. The two coincide when $\theta$ is surjective, in which case
$H^\bullet_{\omni}(\g;\rho)\cong H^\bullet(\g,M)$; this holds for
instance for the adjoint omni-representation
(Corollary~\ref{cor:adjoint-omni-coh}). In general, the two cohomologies
differ. For the first example of Remark~\ref{rmk:M-vs-image}, we have
$H^\bullet_{\omni}(\g;\rho_0)=0$, whereas $H^0(\g,M)=\C$. For the second
one,
\[
\dim H^0(\g,M)=\dim H^0(\g,\g)+1=\dim H^0_{\omni}(\g;\rho_\theta)+1,
\]
these dimensions being finite since $H^0(\g,\g)\subseteq\g/\partial\g$.

(2) Theorem~\ref{thm:graph-omni-cohomology} and
Corollaries~\ref{cor:trivial-omni-coh} and~\ref{cor:adjoint-omni-coh}
identify the coefficient representation $\im(\rho)$ of
Definition~\ref{def:omni-cohomology} with a subrepresentation of the
original module.
\end{remark}

By Proposition~\ref{prop:T-is-ad}, the identity
$T=\ad^L_{M_T}$ identifies
Theorem~\ref{thm:graph-omni-cohomology} with the corresponding
adjoint-type statement for $\g'=M_T$. We retain the formulation in
terms of $T$ because the differential \eqref{eq:graph-delta} is
expressed directly through the embedding tensor.

For a nonzero trivial
omni-representation $\rho_\xi=(0,\xi)$ of $\g$, every omni-cochain
of degree $n\geq1$ is uniquely of the form
\[
\Phi=(0,f),
\qquad
\text{for all}~f\in C^n_{\cLeib}(\g,\C).
\]

\begin{corollary}\label{cor:trivial-omni-coh}
With the above notation,
\[
H^\bullet_{\omni}(\g;\rho_\xi)
\cong H^\bullet(\g),
\]
independently of the choice of the nonzero trivial
omni-representation. If $\xi=0$, then $H^\bullet_{\omni}(\g;\rho_0)=0$.
\end{corollary}

\begin{proof}
If $\xi=0$, then $\im(\rho_0)=0$ and all the cochain spaces vanish.

If $\xi\neq0$, then $\im(\xi)=\C$ and the associated representation
is the trivial representation. Applying
Theorem~\ref{thm:graph-omni-cohomology} with $M=\C$, $T=0$ and
$\theta=\xi$, we obtain the cochain isomorphism
$f\mapsto(0,f)$. All action terms in the differential vanish.
Therefore,
\[
(\delta f)_{\lambda_1,\ldots,\lambda_n}
(x_1,\ldots,x_{n+1})
=
\sum_{1\leq i<j\leq n+1}(-1)^if^{(i,j)},
\]
which is precisely the differential for the trivial representation.
\end{proof}

\begin{remark}\label{rmk:LS-trivial}
The zero omni-representation $\rho_0$ must be treated separately in
Corollary~\ref{cor:trivial-omni-coh}. Indeed, all omni-cochains of
$\rho_0$ vanish, so $H^\bullet_{\omni}(\g;\rho_0)=0$. On the other hand,
$C^0_{\cLeib}(\g,\C)=\C/\partial\C=\C$ and the degree-zero differential
\eqref{eq:delta0} vanishes for the trivial representation, so
$H^0(\g)=\C\neq0$. Hence the isomorphism
$H^\bullet_{\omni}(\g;\rho_\xi)\cong H^\bullet(\g)$ holds only for
$\xi\neq0$.
\end{remark}

\begin{remark}\label{rmk:only-trivial}
Conversely, if $\im(\rho)$ is isomorphic to the trivial representation
$\C$, then $\partial\rho(x)=0$ forces $\varphi=0$, since
$(\partial\varphi(x))_\lambda=-\lambda\varphi(x)_\lambda$, and $\theta(\g)=\pr(\im(\rho))$
is a copy of $\C$ with $\partial=0$. Hence, up to a $\C[\partial]$-linear
embedding $\C\to M$, $\rho$ is a nonzero trivial omni-representation
$\rho_\xi$.
\end{remark}

In the adjoint case, let $\theta:\g\to\g'$ be a homomorphism of
Leibniz conformal algebras, and let
\[
\rho_\theta
=
(\ad^L_{\g'}\circ\theta,\theta).
\]
Then every omni-cochain $\Phi$ of degree $n\geq1$ is uniquely
determined by
$f\in C^n_{\cLeib}(\g,\im(\theta))$ through
\begin{equation}\label{eq:adjoint-cochain}
\Phi=(\ad^L_{\g'}\circ f,f).
\end{equation}
In particular, for
$\Ad=(\ad^L,\Id_{\g})$, every omni-cochain $\Phi$ of degree $n\geq1$
is uniquely of the form
\[
\Phi
=
\Gamma_{\ad^L}\circ f
=
(\ad^L\circ f,f),
\qquad \text{for all }~
f\in C^n_{\cLeib}(\g,\g).
\]

\begin{corollary}\label{cor:adjoint-omni-coh}
With the above notation,
\[
H^\bullet_{\omni}(\g;\rho_\theta)
\cong H^\bullet(\g,\im(\theta)),
\]
where $\im(\theta)$ carries the $\lambda$-actions
\[
x_\lambda v=[\theta(x)_\lambda v]',
\qquad
v_\lambda x=[v_\lambda\theta(x)]'.
\]
In particular,
\[
H^\bullet_{\omni}(\g;\Ad)
\cong H^\bullet(\g,\g).
\]
\end{corollary}

\begin{proof}
By Lemma~\ref{lem:adL-embedding} with $\mathcal A=\g'$,
$\ad^L_{\g'}$ is a conformal embedding tensor and
$\g'_{\ad^L_{\g'}}=\g'$. Applying
Theorem~\ref{thm:graph-omni-cohomology} with
$M=\g'$ and $T=\ad^L_{\g'}$, the map
\[
f\longmapsto(\ad^L_{\g'}\circ f,f)
\]
is an isomorphism of cochain complexes
$C^\bullet_{\cLeib}(\g,\im(\theta))\to C^\bullet_{\omni}(\g;\rho_\theta)$,
which gives the first assertion. For $\g'=\g$ and $\theta=\Id_{\g}$,
we have $\rho_\theta=\Ad$ and $\im(\theta)=\g$ with the adjoint
representation, which gives the second one.
\end{proof}





\end{document}